\documentclass{amsart}
\usepackage{epsf ig}

\newtheorem{theorem}{Theorem}[section]
\newtheorem{lemma}[theorem]{Lemma}
\newtheorem{corollary}[theorem]{Corollary}
\newtheorem{proposition}[theorem]{Proposition}

\theoremstyle{definition}
\newtheorem{definition}[theorem]{Definition}

\theoremstyle{remark}

\numberwithin{equation}{section}

\DeclareMathOperator{\re}{Re}
\DeclareMathOperator{\im}{Im}
\DeclareMathOperator{\supp}{supp}

\DeclareMathOperator{\dv}{div}

\newcommand{\bbC}{{\mathbb C}}
\newcommand{\bbZ}{{\mathbb Z}}
\newcommand{\bbR}{{\mathbb R}}
\newcommand{\bbH}{{\mathbb H}}
\newcommand{\bbN}{{\mathbb N}}
\newcommand{\bbB}{{\mathbb B}}

\newcommand{\calV}{{\mathcal V}}
\newcommand{\calA}{{\mathcal A}}
\newcommand{\calM}{{\mathcal M}}
\newcommand{\calD}{{\mathcal D}}

\newcommand{\calL}{{\mathcal L}}

\newcommand{\calE}{{\mathcal E}}
\newcommand{\calF}{{\mathcal F}}
\newcommand{\calG}{{\mathcal G}}
\newcommand{\calH}{{\mathcal H}}
\newcommand{\calI}{{\mathcal I}}

\newcommand{\del}{\partial}
\newcommand{\bck}{\backslash}
\newcommand{\cinf}{C^\infty}

\newcommand{\tsl}{{\tilde\sigma_l}}
\newcommand{\tsr}{{\tilde\sigma_r}}

\newcommand{\tsigma}{{\tilde\sigma}}
\newcommand{\xx}{X\times X}
\newcommand{\sxx}{{X\mathbin{\times_0}X}}
\newcommand{\esxx}{{X\mathbin{\times_\lambda}X}}
\newcommand{\esxxx}{{(X^3)_\lambda}}

\newcommand{\xU}{\mathbin{\overline\cup}}
\newcommand{\Xl}{X_\lambda}
\newcommand{\bX}{{\del X}}

\newcommand{\xbx}{X\times\bX}

\newcommand{\sxbx}{{X\mathbin{\times_0}\bX}}
\newcommand{\xlxl}{\Xl\times\Xl}
\newcommand{\norm}[1]{\Vert #1\Vert}
\newcommand{\pair}[1]{\langle #1\rangle}

\newcommand{\ds}{d_\tsigma}

\newcommand{\Ds}{D_\tsigma}
\newcommand{\hD}{\hat\Delta}
\newcommand{\ol}{\overline}
\newcommand{\pse}{\psi_\epsilon}
\newcommand{\Dz}{[\Delta - \alpha_0^2\zeta(n-\zeta)]}
\newcommand{\Dl}{\Delta - \lambda}
\newcommand{\hDl}{\hD - \lambda}
\newcommand{\psop}{{^0\Psi}}

\newcommand{\tL}{\tilde\Lambda}

\newcommand{\shlf}{\Omega^{1/2}_0}
\newcommand{\Diff}{\text{\rm Diff}}
\newcommand{\intM}{M^{\circ}}

\begin{document}

\title[Scattering theory for conformally compact metrics]{Scattering theory for conformally
compact metrics with variable curvature at  infinity}
\author{David Borthwick}
\address{Department of Mathematics and Computer Science, 
Emory University, Atlanta}
\email{davidb@mathcs.emory.edu}
\thanks{Supported in part by an NSF Postdoctoral Fellowship.}
\date{October, 2000}



\maketitle
\tableofcontents
\section{Introduction}

Let $X$ be a smooth manifold of dimension $n+1$ with boundary $\bX$, 
equipped with an arbitrary smooth metric
$\bar g$.  A boundary-defining function on $X$ is a function $x \ge 0$ such that $\bX =
\{x=0\}$ and $dx\ne 0$ on $\bX$.  A \textit{conformally compact} metric on the interior of $X$
is a metric of the form
$$
g = \frac{\bar g}{x^2}.
$$
Such metrics were introduced by Mazzeo \cite{Ma88} as a generalization of the hyperbolic metric
on $\bbB^n$.  The metric $g$ is necessarily complete.  Any non-trapped geodesic $\gamma$
approaches a point $y\in \bX$, and as $t\to\infty$ all sectional curvatures at $\gamma(t)$ all
approach the value $-\alpha(y)^2$, where
$$
\alpha = |dx|_{\bar g} \text{ restricted to }\bX.
$$
Let $\alpha_0 = \inf_{\bX}\alpha$ and $\alpha_1 = \sup_{\bX}\alpha$.

Mazzeo established the basic properties of the spectrum of the
$p$-form Laplacian associated to $g$ and proved the appropriate Hodge theorem for this
context.  Here we will deal only with the Laplacian on functions, denoted simply
by $\Delta$ since $g$ is fixed throughout.

\begin{theorem}\label{mathm}\cite{Ma88, Ma91a}
The essential spectrum of $\Delta$ is $[\frac{\alpha_0^2n^2}4, \infty)$ and is absolutely
continuous. There are no embedded eigenvalues except possibly at $\frac{\alpha_0^2n^2}4$.
\end{theorem}

The case where $\alpha = \alpha_0$ is constant is referred to as \textit{asymptotically
hyperbolic}, since sectional curvatures all approach $-\alpha_0^2$ at infinity.  Note that
constant curvature ``at infinity'' does not imply that $g$ has constant curvature at
any point.

For asymptotically hyperbolic metrics, Mazzeo-Melrose proved meromorphic continuation of the
resolvent.
\begin{theorem}\label{mmthm}\cite{MM}
If $\alpha = \alpha_0$ then the resolvent $R_\zeta = [\Delta - \alpha_0^2\zeta(n-\zeta)]^{-1}$
has a meromorphic continuation to $\bbC \bck \frac12(n - \bbN)$.
\end{theorem}

The proof is by a parametrix construction which gives a detailed picture of the structure of
the resolvent kernel.  One sees, for example, that 
$$
R_\zeta: \dot\cinf(X) \to x^{\zeta}\cinf(X),
$$ 
where $\dot\cinf(X)$ is the space of smooth functions vanishing to
infinite order at $\bX$.  This property leads to a ``functional parametrization'' of the
continuous spectrum (see \cite{MeGS}).  Given
$f|_\bX\in\cinf(\bX)$, we can solve away a Taylor series at the boundary to extend
$f$ smoothly into the interior in such a way that $\Dz x^{n-\zeta} f \in \dot\cinf(X)$.  Then by
applying $R_\zeta$ to the remainder we construct a generalized eigenfunction $u$ solving
$\Dz u = 0$ with 
\begin{equation}\label{uff}
u = x^{n-\zeta} f + x^\zeta f',
\end{equation}
where $f'\in\cinf(X)$.  This $u$ is uniquely determined by $f|_\bX$, and the map $E_\zeta: f|_\bX \mapsto u$ is
called the Poisson operator, after the classical case.   It defines a parametrization of the continuous
spectrum by $\cinf(\bX)$.  This construction also yields the scattering matrix 
$S_\zeta:f|_\bX \mapsto f'|_\bX$, which is a pseudodifferential operator of order $2\re\zeta-n$.
Note that as defined here, $E_\zeta$ and $S_\zeta$ depend on the choice of $x$.  This dependence is easily
removed by considering sections of a certain trivial line bundle instead of functions, but for out
purposes it is much more convenient to fix a choice of $x$ for the whole paper.
The kernels of $E_\zeta$ and $S_\zeta$ can be derived from $R_\zeta$ and are meromorphic functions of
$\zeta \in \bbC\bck\frac12(n - \bbN)$.

This paper is devoted to the extension of such results to the general conformally compact case, with variable
$\alpha$.   To heighten the analogy with the asymptotically hyperbolic case, we continue to use a 
spectral parameter $\zeta$ such that the relation to the eigenvalue $\lambda$ is
$$
\lambda = \alpha_0^2\zeta(n-\zeta).
$$ 
For constant $\alpha$ this association comes from the equation for indicial roots of $\Dl$ (see
\S\ref{indsec}).  In the general case the indicial roots are variable, and even singular for certain 
values of $\lambda$.   This complication is the source of interesting new features in the scattering 
theory.

\begin{figure}
\centerline{\epsfysize=1.5in \epsfbox{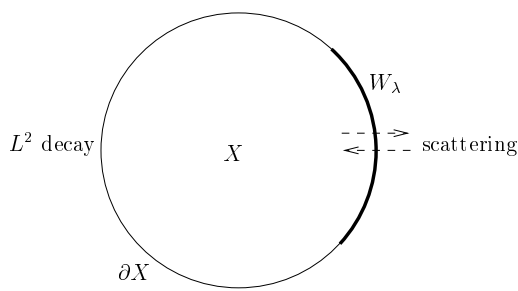}}
\caption{Behavior of generalized eigenfunctions
for the irregular continuous spectrum. Scattering occurs
only where the curvature at infinity is $> -\frac{4\lambda}{n^2}$.}
\label{scdiagram}
\end{figure}

To summarize the results of this paper:
\begin{enumerate}
\item  The Mazzeo-Melrose parametrix construction can be used to obtain meromorphic continuation of the 
resolvent $\Dz^{-1}$ to the plane minus a set $\Gamma$ which is a collection of intervals (Theorem
\ref{merores}).   Figure \ref{mero} shows the region of meromorphic continuation, which
is defined by the condition that the indicial root avoid the set
$\frac12(n-\bbN_0)$ and includes the portion of the continuous spectrum $\lambda \in (\frac{\alpha_1^2n^2}4,
\infty)$.  Within the region of meromorphic continuation we can construct the Poisson kernel and scattering
matrix as in  the asymptotically hyperbolic case (Propositions \ref{poissk} and \ref{smatrix}). 
\begin{figure}
\centerline{\epsfysize=1.7in \epsfbox{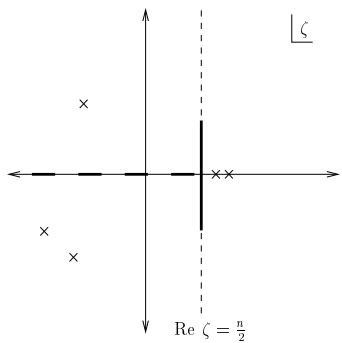}}
\caption{The resolvent is meromorphic outside the marked intervals.  The dotted line indicates
the continuous spectrum, and $\times$'s represent possible poles.}
\label{mero}
\end{figure}
\item  We'll refer to the portion of the continuous spectrum inaccessible by meromorphic continuation,
$\lambda\in [\frac{\alpha_0^2n^2}4, \frac{\alpha_1^2n^2}4]$ as the \textit{irregular} continuous spectrum. 
(This corresponds to the vertical solid line in Figure \ref{mero}.)
We analyze it through a limiting absorption principle.  This means establishing the strong limit of
$\Dz^{-1}$ as $\zeta$ approaches the line $\re\zeta = \frac{n}2$ from the right (Theorem
\ref{lathm}).   With this technique we show that the \textit{scattering set}
$$
W_\lambda = \{\alpha^2 < \tfrac{4\lambda}{n^2}\} \subset \bX,
$$
``parametrizes'' the irregular continuous spectrum at $\lambda$.  More precisely, for
$f\in\dot\cinf(W_\lambda)$ there is a unique solution of $\Dz u = 0$ having an asymptotic expansion with
leading behavior
$$
u \sim x^{n-\sigma(\zeta,y)} f(y) +  x^{\sigma(\zeta,y)} f'(y) \text{ near }W_\lambda,
$$
where $\sigma$ is the (variable) indicial root defined in \S\ref{indsec} and $f'\in \cinf(W_\lambda)$
(Proposition \ref{nuparam}).  The restriction of $u$ to a neighborhood intersecting the boundary only in
$\{\alpha^2 > \tfrac{4\lambda}{n^2}\}$  is in $L^2$.  
Thus, as illustrated in Figure \ref{scdiagram},  generalized eigenfunctions only ``propagate out to
infinity'' on $W_\lambda$.    By combining the local parametrix construction with the
limiting absorption principle, we are able to construct the Poisson kernel, understand its structure near
$W_\lambda$, and thus show that the scattering matrix is a pseudodifferential operator defined on
$W_\lambda$ (Theorem \ref{nusmatrix}).
\item  The edge of scattering set is the \textit{crossover region} 
$\{\alpha^2 = \frac{4\lambda}{n^2}\} \subset\bX$.  
If we assume that $\frac{2\sqrt{\lambda}}n$ is a regular value of $\alpha$, 
so that the crossover region is a submanifold of $\bX$, then we can undertake a direct construction
of the resolvent $\Dz^{-1}$.   
The technique is to blow up the crossover region in $X$ to resolve the singularities of the indicial
root.   Adapting the parametrix construction to this extra blow-up, we obtain a good parametrix
from which the structure of the resolvent may be deduced.
The result is a full picture of the boundary behavior of
the resolvent kernel under this generic assumption on $\lambda$ (Theorem \ref{resdirect}).
\end{enumerate}

This behavior of generalized eigenfunctions can be interpreted physically.
Larger $\alpha$ corresponds to more rapid volume growth at infinity, so one would naturally expect waves 
traveling in such directions to diffuse more quickly.  
At low frequencies (relative to $\alpha$)
the diffusion effect is evidently strong enough to overcome propagation, while sufficiently
high-frequency waves do propagate in all directions.

Scattering theory on hyperbolic manifolds has an extensive literature (see \cite{Hi} for a review 
of the subject).   For this case the absence of embedded eigenvalues was proven by Lax-Phillips in
\cite{LP82} and meromorphic continuation of the resolvent by Perry \cite{P89}, independently of \cite{MM}. 
Perry also proved that the scattering operator was pseudodifferential and computed its symbol.

For asymptotically hyperbolic metrics, the parametrization of the spectrum as in (\ref{uff}),
and the corresponding definition of the scattering matrix, again a pseudodifferential operator, 
was implicit in \cite{MM} (stated, for example, in \cite{MeGS}).  The proof was given by 
Joshi-S\'a Barreto \cite{JS}, who prove an inverse result on the determination of asymptotics
of the metric from the symbol of the scattering operator.  The equivalence of resolvent and scattering
resonances for asymptotically hyperbolic metrics was proven by Borthwick-Perry \cite{BP}.  In \cite{GZ95}
Guillop\'e-Zworski establish an upper bound on the counting function for resonances, under the
stronger assumption of hyperbolic `near infinity,' i.e. outside a compact set.  For the case of Einstein
metrics which are asymptotically hyperbolic Lee proved that there is no discrete spectrum
provided the Yamabe invariant of the induced conformal structure on $\bX$ is non-negative \cite{L95}.

Outside of \cite{Ma88, Ma91a} no work seems to have been done on conformally compact metrics in full
generality.  The phenomenon of a scattering operator defined for a frequency-dependent set
of directions appears to be quite new, although the direction-dependence bears some analogy with 
recent results of Herbst-Skibsted on scattering by homogeneous potentials \cite{HS}.

\bigskip\noindent{\bf Acknowledgments.}  Thanks to Rafe Mazzeo for encouragement
to undertake this project, and Richard Melrose for some very helpful suggestions.
The work was stimulated by a workshop at the Stefan Banach International Mathematical Centre in
Warsaw, August 1999, for which support is gratefully acknowledged.

\section{Boundary asymptotics}\label{basymsec}

Our main goal is to describe precisely the structure of the generalized functions, as well as the 
resolvent, Poisson, and scattering kernels.  The first step is to introduce the spaces
which will characterize the behavior of these functions at the boundary.

Let $M$ be a smooth manifold with corners (see \cite{Me92, MeDA} for basic definitions). 
The boundary hypersurfaces of $M$, themselves manifolds with corners, are labeled
$Y_j$, $j=1,\dots,p$, and we introduce corresponding boundary defining functions $\rho_1,\dots,\rho_p$.  
We are basically interested in functions which behave near $Y_j$ 
like $\rho_j^{\beta_j}$ times a smooth function, for a set of index functions $\beta_j \in \cinf(Y_j)$.   
However, to make a definition independent of the choice of $\rho_j$'s, `smooth' must be
relaxed to `polyhomogeneous conormal' with a particular index set.  We'll follow closely the
definition of spaces of polyhomogeneous conormal functions in
\cite{Me92, MeDA}, but give a self-contained presentation for the convenience of the reader.

The set of smooth vector fields tangent to the boundary is denoted by $\calV_b(M)$.  
As an auxiliary space in the definition, define for a multi-index $m\in \bbR^p$ the space
$$
\calA^q(M) = \{u\in \cinf(\intM): (\calV_b)^k u \in \rho^m L^\infty(M)\; \forall k\},
$$
where $\intM$ denotes the interior and $\rho^m = \rho_1^{m_1} \dots \rho_p^{m_p}$.
This is clearly invariant under the action of $\calV_b$.
Because of logarithmic terms it will be convenient to use the space
$$
\calA^{m-} = \bigcap_{m'<m} \calA^{m'},
$$
where $m'< m$ means $m'_j<m_j$ for each $j$.

Given $\rho_j$ we choose a product decomposition 
$Y_j \times [0,\epsilon)_{\rho_j}$ of a neighborhood of $Y_j$ in $X$. 
Within this product neighborhood the radial vector field is
$$
V_j = \rho_j \del_{\rho_j},
$$
which we'll extend to the rest of $M$ so as to define an element of $\calV_b(M)$.
$V_j$ is determined independently of the product decomposition up to an element of $\rho_j \calV_b(M)$.
\begin{definition}\label{phgvo}
For a family of smooth functions $\beta\in \cinf(M;\bbR^p)$, the space
$\calA_\beta(M)$ of `polyhomogeneous functions with variable order' consists of functions 
$u\in\cinf(\intM)$ such that for any $m\in \bbN_0^p$ we have
$$
\Bigl[\prod_{l=1}^p \prod_{k=0}^{m_l-1} (V_j - k)^{k+1}\Bigr] (\rho^{-\beta}u)
\in \calA^{m-}(M).
$$
\end{definition}

To signify vanishing to infinite order at a particular boundary face we'll use the notation $\beta_j=\infty$.
The invariance of $\calA_\beta(M)$ under the action of $\calV_b(M)$ is immediate from the definition.  
When $\beta$ is constant it is usual to define the spaces with operators $(V_j - \beta_j - k)$
acting on $u$.  This is not equivalent here, unless each $\beta_j$ is independent of $\rho_j$
near $\rho_j=0$.  We will see below that one could always make this assumption, because $\calA_\beta(M)$
depends only on $\beta_j|_{Y_j}$.

The space $\calA_\beta(M)$ could just as well be characterized by the existence of asymptotic expansions.
Near the boundary surface $Y_j$ the expansion will take the form
\begin{equation}
u \sim \sum_{0\le l\le k<\infty} \rho_j^{\beta_j+k} (\log\rho_j)^l a_{k,l},
\end{equation}
for functions $a_{k,l}$ on $Y_j$.  By this we mean that for any $q>0$ 
\begin{equation}\label{aedef}
\rho^{-\beta}u - \sum_{0\le l\le k<q-1} \rho_j^{k} (\log\rho_j)^l \chi(\rho_j) a_{k,l}
\in \calA^{q\pair{j}-}(M),
\end{equation}
where $q\pair{j}$ denotes the index set $(\dots,0,q,0,\dots)$ with $q$ in the $j$-th place,
and $\chi\in \cinf([0,\infty))$ with $\chi = 1$ on $[0,\epsilon/2]$ and $\chi=0$ on $[\epsilon,\infty)$,
so that $ \chi(\rho_j) a_{k,l}$ may be thought of as a function on $M$ which vanishes outside the product
neighborhood

To a multi-index $\beta$ on $M$ we can associate a multi-index $\beta^{(j)}$ on each face $Y_j$.  
If $H_l$ is a boundary hypersurface of $Y_j$ (and hence a corner of $M$), then set
$\beta^{(j)}_l = \beta_k|_{H_l}$ where $Y_k$ is the unique boundary surface such that $H_l$ is a
component of $Y_j\cap Y_k$.  

\begin{proposition}\label{bexp}
If $u\in \cinf(\intM)$, then $u\in \calA_\beta(M)$ if and only if $u$ has an asymptotic expansion
at each boundary surface $Y_j$:
$$
u \sim \sum_{0\le l\le k<\infty} \rho_j^{\beta_j+k} (\log\rho_j)^l a_{k,l},
$$
where $a_{k,l} \in \calA_{\beta^{(j)}}(Y_j)$
\end{proposition}
\begin{proof}
It suffices to prove the expansion for $\rho^{-\beta} u$, so we can assume $\beta=0$.  This
is a special case of a result proven in \cite{MeDA} using the Mellin transform.  
We'll give a different proof using methods found in \cite{Jo97, JS}.

It also suffices to consider a single face, say $Y_1$.  We'll work in a product neighborhood $Y_1\times
[0,1)_{t}$ and ignore the cutoff $\chi$ in (\ref{aedef}).
From the definition we have
\begin{equation}\label{startpt}
\Bigl[\prod_{k=0}^{q} (t\del_t - k)^{k+1}\Bigr] u
\in \calA^{(q+1)\pair1-}(M).
\end{equation}
Let 
\begin{equation}\label{u1def}
u_1 = \Bigl[\prod_{k=1}^q (t\del_t - k)^{k+1} \Bigr] u,
\end{equation}
so that (\ref{startpt}) becomes the estimate 
\begin{equation}\label{dtau}
\del_t(u_1) = f \in \calA^{q\pair1-}(M).
\end{equation}
In particular, since $q>0$ we see that $u_1$ approaches a limit as $t\to 0$, call it $b_1$.
This $b_1$ is a function on $Y_1$, and since we also have estimates of the form
$$
\del_t(\calV_b)^k u_1 \in \calA^{q\pair1-}(M),
$$
we have uniform convergence of tangential derivatives $(\calV_b)^k u_1$ as $t\to 0$, giving us
bounds on $(\calV_b)^k b_1$.   This shows that $b_1\in \calA^0(Y_1)$.  
The same principle applies to the radial vector field estimates.  
For example, since
$$
V_2 u_1 \in \calA^{1\pair2-}(M) \text{ and }\del_t V_2 u_1 \in \calA^{(q,1,0,\dots)-}(M)
$$
the limit of $V_2 u_1$ exists as $t\to 0$ and equals $V_2|_{Y_1} b_1$.  In this way
we get an estimate $V_2|_{Y_1} b_1 < C(\rho_2|_{Y_1})^{1-\delta}$ for any $\delta>0$.  
Repeating this argument with higher derivatives and at all boundary faces of $Y_1$, we conclude that
$$
b_1 \in \calA_{0}(Y_1).
$$
Integrating (\ref{dtau}) from $0$ to $t$ gives
$$
u_1 - b_1 = \int_0^t f\>dt.
$$
Again using $q>0$, it is easy to see that
$$
\int_0^t f\>dt \in \calA^{(q+1)\pair1-}(M),
$$
and hence $u_1 - b_1\in \calA^{(q+1)\pair1-}(M)$.
Setting $a_{0,0} = b_1/\prod_{k=1}^{q} (-k)^{k+1}$, we thus derive from (\ref{u1def}) that
\begin{equation}\label{k1est}
\Bigl[\prod_{k=1}^{q} (t\del_t -k)^{k+1} \Bigr] (u - a_{0,0}) \in \calA^{(q+1)\pair1-}(M)
\end{equation}

Now let
$$
u_2 = \Bigl[\prod_{k=2}^q (t\del_t - k)^{k+1} \Bigr] (u - a_{0,0}),
$$
so that (\ref{k1est}) implies
$$
(t\del_t-1)^2 u_2 \in \calA^{(q+1)\pair1-}(M).
$$
This could be written 
$$
\del_t [t^{-1} (t\del_t-1) u_2] \in \calA^{(q-1)\pair1-}(M),
$$
so assuming $q>1$ we apply the above argument to find $b_2\in\calA^{0}(Y_1)$ such that
$$
t^{-1} (t\del_t - 1)u_2 - b_2 \in \calA^{q\pair1-}(M).
$$
Rewrite this as
$$
\del_t [t^{-1} u_2 -  (\log t) b_2] \in \calA^{(q-1)\pair1-}(M).
$$
Repeating the argument yet again, we find $b_2'\in\calA^{0}(Y_1)$ such that
$$
u_2 - (t\log t) b_2 - t b_2' \in \calA^{(q+1)\pair1-}(M).
$$
Then from $b_2$ and $b_2'$ we form $a_{1,1}$ and $a_{1,0}$ by linear combination so that
$$
\Bigl[\prod_{k=2}^{q} (t\del_t -k)^{k+1} \Bigr] \Bigl(u - a_{0,0} - (t\log t) a_{1,1} - t a_{1,0}\Bigr) 
\in \calA^{(q+1)\pair1-}(M)
$$

This procedure may be continued inductively up to $k=q-1$, yielding
\begin{equation}\label{ttasym}
(t\del_t - q)^{q+1} \Bigl(u - \sum_{j=0}^{q-1} \sum_{l=0}^j 
t^{j} (\log t)^l a_{j,l} \Bigr) \in \calA^{(q+1)\pair1-}(M).
\end{equation}
To remove the remaining derivatives, consider a function $F\in\cinf(M^\circ)$ satisfying
\begin{equation}\label{Fstart}
(t\del_t - q)F \in \calA^{q\pair1-}(M),
\end{equation}
which is equivalent to
$$
\del_t(t^{-q}F) = h\in \calA^{(-1)\pair1-}(M).
$$
Integrating from $t$ to $1$ gives
$$
t^{-q}F = f(1) - \int_t^1 h\>dt.
$$
Since we have $(\calV_b)^k \int_t^1 h\>dt = O(\log t)$ by the estimates on $h$,
we conclude that (\ref{Fstart}) implies $F\in \calA^{q\pair1-}(M)$.

Applying this argument $(q+1)$ times to (\ref{ttasym}) we obtain
$$
u - \sum_{j=0}^{q-1} \sum_{l=0}^j 
t^{j} (\log t)^l a_{j,l}  \in \calA^{q\pair1-}(M),
$$
where $a_{j,l} \in \calA_0(Y_1)$.
\end{proof}

We will need a means to indicate subspaces of $\calA_\beta$ consisting of functions with truncated expansions,
i.e. with a certain number of the leading coefficients set to zero.
For $q\in \bbN^p$ define
\begin{equation}\label{truncdef}
\calA_{\beta| q}(M) = \prod_{l=1}^p (\rho_l\ln\rho_l)^{q_l} \cdot \calA_\beta(M). 
\end{equation}
A useful alternative characterization is
\begin{equation}\label{truncdef2}
\calA_{\beta| q}(M) = \rho^\beta 
\Bigl[\prod_{l=1}^p \prod_{k=0}^{q_l-1} (V_j - k)^{k+1}\Bigr] \rho^{-\beta} \cdot \calA_\beta(M),
\end{equation}
The proof of Proposition \ref{bexp} shows that (\ref{truncdef}) and (\ref{truncdef2}) are equivalent.

Generally one defines spaces of polyhomogeneous conormal functions with more general expansions,
specifying explicitly the set of possible powers of $\rho_j$ and $\log\rho_j$ at each face.  One could do
the same for variable orders, but the possible crossing of orders complicates the definition.
Since we do not require such generality, in Definition \ref{phgvo} we have chosen to use the minimal set of
powers consistent with variable order.
\begin{lemma}\label{aindep}
The space $\calA_\beta(M)$ is independent of the choice of $V_j$ and depends on $\beta$ only through the
restrictions $\beta_j|_{Y_j}$.
\end{lemma}
\begin{proof}
Since $V_j$ is unique up to  $\rho_j \calV_b(M)$, independence follows from the invariance of $\calA^q$
under $\calV_b$.

To study the $\beta$ dependence, for simplicity let us specialize to the case $p=1$, i.e. a manifold
with boundary $X$ with boundary defining function $x$.  It suffices to consider an index $\beta = xf$ where
$f\in\cinf(X)$, and show that
$\calA_{0}(X) \subset \calA_{xf}(X)$.  Define the commutator operators
$$
Z_q = \Bigl[\prod_{k=0}^{q-1} (x\del_x - k)^{k+1}, x^{(xf)}\Bigr].
$$
By induction we will show that
\begin{equation}\label{zqmap}
Z_q: \calA_0(X) \to \calA_{0|q}(X).
\end{equation}
Because $\calA_{0|q}(X) \subset \calA^{q-}(X)$, for $u\in\calA_0(X)$ we can then estimate 
\begin{equation*}\begin{split}
\prod_{k=0}^{m-1} (x\del_x - k)^{k+1}(x^{(xf)}u) &= Z_m u + x^{(xf)} \prod_{k=0}^{m-1} 
(x\del_x - k)^{k+1} u \\
&\in \calA^{m-}(X),
\end{split}\end{equation*}
and hence $u\in\calA_{xf}$.

So the result follows once we establish (\ref{zqmap}).  $Z_1$ is multiplicative,
$$
Z_1 = x^{(xf)} [(x\log x+x)f + x^2\log x (\del_x f)],
$$
and the mapping property $Z_1:\calA_0(X)\to \calA_{0|1}(X)$ is immediate.
For the inductive step, assume $Z_q$ satisfies (\ref{zqmap}) and consider
$$
Z_{q+1} u = (x\del_x - q)^{q+1} Z_{q}u + [(x\del_x-q)^{q+1}, x^{(xf)}] \prod_{k=0}^{q-1} (x\del_x - k)^{k+1} u
$$
By definition we have
$$
(x\del_x - q)^{q+1}: \calA_{0| q}(X) \mapsto \calA_{0|(q+1)}(X), 
$$
which takes care of the first term.  And the second is a sum of terms of the form
$$
(\calV_b)^k\cdot Z_1\cdot (\calV_b)^{q-k} \prod_{k=0}^{q-1} (x\del_x - k)^{k+1} u,
$$
and so using (\ref{truncdef2}) the desired mapping follows from $Z_1:\calA_{0|q}(X) \to \calA_{0|(q+1)}(X)$.

The argument is the same in the general case, except that the induction must be
done over $\bbN^p$.
\end{proof}

\section{Indicial roots}\label{indsec}

It will be convenient to put the metric into a normal form.  The proof of the following result
may be adapted directly from \cite{JS}, where $\alpha=1$.
\begin{proposition}
Let $X$ be a compact manifold with $g$ a conformally compact metric.  There exists a product
decomposition $(x,y)$ near $\bX$ such that
$$
g = \frac{dx^2}{\alpha^2(y) x^2} + \frac{h(x,y,dy)}{x^2} + O(x^\infty).
$$
Here $-\alpha(y)^2$ is the limiting curvature at infinity.
\end{proposition}

For the rest of the paper we will assume that 
\begin{equation}\label{norform}
g = \frac{dx^2}{\alpha^2(y) x^2} + \frac{h(x,y,dy)}{x^2},
\end{equation}
in some product neighborhood of the boundary.
The $O(x^\infty)$ correction is dropped for notational simplicity, since its effect
on the various estimates and asymptotic expansions at the boundary would be trivial.

The corresponding Laplacian operator is
\begin{equation}
\Delta = \alpha^2 \Bigl[-(x\del_x)^2 + nx\del_x - x^2(\del_x \log
\sqrt{h}) \del_x\Bigr] + x^2\Delta_h -x^2 (\del_i \log\alpha) h^{ij} \del_j.
\end{equation}
The \textit{indicial equation} in this context is obtained by setting the leading order term in
$(\Dl)x^\sigma$ equal to zero and solving for the indicial root $\sigma$ as a function of $\lambda$.
This yields
$$
\alpha^2 \sigma(n-\sigma) - \lambda = 0,
$$
so $\sigma$ depends on $y$ through $\alpha(y)$.

\begin{figure}
\centerline{\epsfysize=1.7in \epsfbox{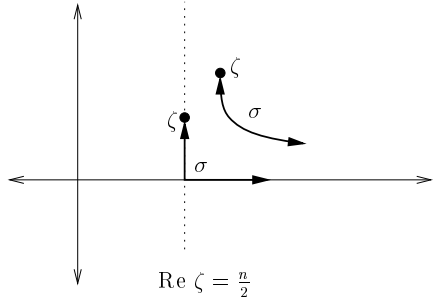}}
\caption{For fixed $\zeta$, the indicial root $\sigma$ varies with $y\in\bX$.}
\label{indroot}
\end{figure}

In the hyperbolic case it is natural and customary to choose spectral parameter $\zeta$ equal to
the indicial root.  For variable $\alpha$ it seems natural to use the relation  
\begin{equation}
\lambda = \alpha_0^2 \zeta(n-\zeta),
\end{equation}
so that at least the continuous spectrum still corresponds to the line $\re\zeta=\frac{n}2$.
Then the indicial equation can be solved for $\sigma(\zeta,y)$ 
\begin{equation}
\sigma = \sigma(\zeta, y) = \frac{n}2 + \sqrt{\frac{n^2}4 - \frac{\alpha_0^2}{\alpha^2}
\zeta(n-\zeta)}.
\end{equation}
The square root has a natural interpretation such that $\re\sigma>n/2$ whenever $\re\zeta>n/2$,
which breaks down only if $\frac{n^2}4 - \frac{\alpha_0^2}{\alpha^2}
\zeta(n-\zeta)=0$ for some $y\in\bX$.  We can thereby define $\sigma$ as an analytic function
of $\zeta$ outside of the segment  
$$
\zeta \notin\Bigl\{ \zeta=\frac{n}2+it:\;  t^2
\le\frac{n^2}4 \bigl(\frac{\alpha_1^2}{\alpha_0^2}-1\bigr)\Bigr\}.
$$  
The range of definition is extended to $\zeta\in\bbC$ by continuity, and
Figure \ref{indroot} shows the behavior of $\sigma(\zeta)$
as a function of $y$.  The indicial root fails to be analytic at those values of $\zeta$ for which
$\lambda \in [\frac{\alpha_0^2n^2}4, \frac{\alpha_1^2n^2}4]$.

It turns out that meromorphic continuation of the resolvent requires not only that $\sigma$ be an analytic
function of $\zeta$, but that $\sigma$ avoid the points $\frac12(n-\bbN)$.  Accordingly we define
$$
\Gamma = \{\zeta\in\bbC:\; \sigma(\zeta,y) \in \tfrac12(n-\bbN_0) \text{ for some }y\in\bX\},
$$
as pictured in Figure \ref{mero}.
Since $\zeta\neq \frac{n}2$ implies $\sigma(\zeta,\cdot)\in\cinf(\bX)$,
we may define $\calA_\sigma(X)$ for $\zeta\in\bbC\bck\Gamma$ by choosing an arbitrary smooth extension
of $\sigma(\zeta,\cdot)$ off the boundary.   (We will consider the case of singular $\sigma$ starting in
\S\ref{lasec}.)  

\begin{lemma}\label{solveaway}
Let $\zeta\in\bbC\bck\Gamma$.
Given $v\in\calA_{\sigma|1}$, we can find $u\in\calA_{\sigma|1}$ such that
$$
v - \Dz u \in\dot\cinf(X).
$$
\end{lemma}

\begin{proof}
The leading terms in the boundary expansion of $v$ are
$$
x^{\sigma+1} (\log x) g_1 + x^{\sigma+1} g_2,
$$  
where $g_1,g_2\in \cinf(X)$.
Observe that for any $\psi\in\cinf(X)$
$$
\Dz x^{\sigma+1} (\log x) g_1 = \alpha^2(n - 2\sigma -1) x^{\sigma+1} (\log x) g_1 + x^{\sigma+1} h_2
+ R,
$$
where $h_2\in\cinf(X)$ and $R\in \calA_{\sigma|2}$.
Provided that $\sigma \ne (n-1)/2$ we can set
$$
u_1 = \frac{1}{\alpha^2(n-2\sigma-1)} \Bigl[x^{\sigma+1} (\log x) g_1 + x^{\sigma+1} (g_2-h_2)\Bigr],
$$
so as to have
$$
v - \Dz u_1 \in \calA_{\sigma|2}
$$

The remaining terms $u_j$ are obtained by an obvious induction, with the requirement that $\sigma \notin
\frac12(n-\bbN_0)$ ensuring that no zeroes occur in denominators. 
Then using Borel's lemma we sum the series asymptotically at $x=0$ to get $u\sim \sum u_j$.
\end{proof}

Since $\sigma$ is the indicial root, $\Dz x^\sigma f \in \calA_{\sigma|1}$
for $f\in\cinf(X)$, and we immediately conclude the following:
\begin{corollary}\label{cdsolve}
Let $\zeta\in\bbC\bck\Gamma$.
Given $f \in \cinf(\bX)$  we can solve
$$
\Dz u \in \dot\cinf(X),
$$
for $u\in\calA_\sigma(X)$ such $u$ has leading boundary term 
$$
u(x,y) \sim x^{\sigma} f(y) \text{ as } x\to 0.
$$
\end{corollary}

\section{Parametrix construction}\label{parsec}

The operator $\Delta$ is a member of a class of differential
operators $\Diff_0^m(X)$ (where $m$ denotes the order) generated by $\calV_0(X)$, 
the set of smooth vector fields on $X$ which vanish at the boundary.  
In local coordinates $\calV_0$ is generated by $x\del_x$ and $x\del_y$ (whereas
$\calV_b$ is generated by $x\del_x$ and $\del_y$).

\begin{figure}
\centerline{\epsfysize=2.3in \epsfbox{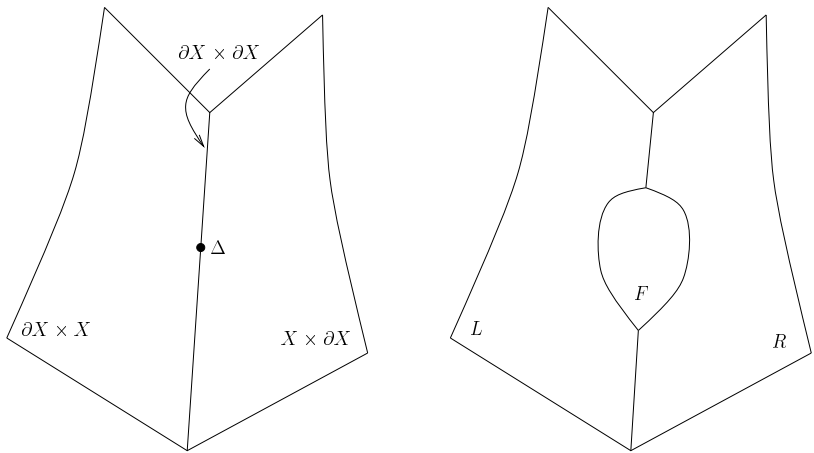}}
\caption{The diagonal $S$ in the corner of $\xx$ is blown up to form
the stretched product $\sxx$.}
\label{sxx}
\end{figure}

The \textit{stretched product} $\sxx$ was introduced in \cite{MM, Ma88} as the natural space on
which to study integral kernels of operators in $\Diff_0^m(X)$.  To define it,
let $S = \Delta(\bX\times\bX)\subset \xx$, which is the intersection of the diagonal with the corner.   
The stretched product is formed by blowing up this submanifold, which is notated:
$$
\sxx = [\xx; S].
$$ 
As a set, $\sxx$ is $\xx$ with $S$ replaced by a \textit{front
face} consisting of its (inward-pointing) spherical normal bundle, a procedure best thought of as the
introduction of polar coordinates around $S$.   The blow-up is illustrated in Figure \ref{sxx}.
In local coordinates $(x,y,x',y')$ for $\xx$ corresponding to the product decomposition of $X$
near $\bX$ the diagonal in the corner is $S = \{x=x'=y-y'=0\}$.
Letting $r = \sqrt{x^2+{x'}^2 + (y-y')^2}$, we introduce the coordinates $(\rho,\rho',\omega,r,y)$ for
$\sxx$, where 
$$
\rho = \frac{x}r,\qquad \rho' = \frac{x'}r,\qquad \omega = \frac{y-y'}r.
$$
The boundary faces of $\sxx$ are the front face $r=0$, the \textit{left face} $\rho=0$,
and the \textit{right face} $\rho'=0$.  The blow-down map is denoted by 
$$
\beta : \sxx\to\xx.
$$

The operators occuring in the parametrix construction are
characterized by the behavior of the lifts of their distribution kernels to $\sxx$.  
Since $g$ is fixed, we can associate to each
operator its distributional kernel with respect to the Riemannian density.
(Using half-densities avoids this dependency and is thus better for many purposes, but it would introduce
another layer of notation which we don't actually need at this point.)  

The microlocalization of $\Diff_0^m(X)$ is $\psop^m(X)$, the space of pseudodifferential operators
modeled on $\calV_0(X)$.  This was introduced in \cite{MM} and is often referred to as the `small
calculus.'   An operator is in $\psop^m(X)$ if the lift of its kernel to $\sxx$ has
conormal singularity of order $m$ at the lifted diagonal which is extendible over the double across the front
face.  The lifted kernels are further required to vanish to infinite order at the left and right faces.

Let $\sigma_l$ and $\sigma_r$ be the pullbacks of $\sigma$ through $\bX\times X$ and $X\times\bX$ and
up to the left and right faces, respectively.  We order the faces of $\sxx$ left, right, and front.
Under the assumption that $\sigma$ is smooth we introduce the space
$\psop_{\sigma_l,\sigma_r}(X)$ of operators whose kernels satisfy
$$
\beta^*K \in \calA_{\sigma_l,\sigma_r,0}(\sxx),
$$
but which furthermore are extendible across the front face (hence no logarithmic terms in the expansion
there). The residual class in the construction is $\Psi_{\sigma_l, \sigma_r}(X)$, consisting of operators
with kernels in $\calA_{\sigma_l,\sigma_r}(\xx)$.

For future reference, we record the following mapping properties.
\begin{lemma}\label{maplemma}
\begin{equation*}\begin{split}
\psop^{m}(X)&:\dot\cinf(X) \to \dot\cinf(X) \\
\psop_{\sigma_l,\sigma_r}(X)&: \dot\cinf(X)\to \calA_\sigma(X)
\end{split}\end{equation*} 
\end{lemma}
\begin{proof}
The first property follows from a standard wave-front set argument.  To prove the second,
consider $A\in \psop_{\sigma_l,\sigma_r}(X): \dot\cinf(X)$ and $f\in \dot\cinf(X)$.
We can compute $Af$ by first pull $f$ up to $\sxx$ through the right, then multiplying by the lift of $A$ 
times the Riemannian density (in the right factor), then pushing forward to $X$ through the left.  But the
lift of $f$ vanishes to infinite order at both the right and front faces, and so the
push-forward can in fact be written as an integral on $\xx$:
$$
Af(x,y) = \int F(x,y, x',y') \>dg(x',y'),
$$
for some $F \in \calA_{\sigma_l,\infty}(X\times X)$.  Then $Af \in \calA_\sigma(X)$
is established by moving derivatives under the integral.
\end{proof}

As noted in \cite{Ma91a}, the Mazzeo-Melrose parametrix construction of \cite{MM} applies
locally to the case of variable $\alpha$ without much alteration, provided we restrict to 
a neighborhood of the boundary where $\sigma\notin \frac12(n-\bbN_0)$.  
The only real change in the construction is the addition of logarithmic terms.  
Lemma \ref{solveaway} shows shows that these are easily handled when solving away Taylor 
series at the boundary.  In addition, we need to include logarithmic singularities in
applications of the model hyperbolic resolvent on fibers of the front face, extending Proposition 6.19 of
\cite{MM}.   

The model case in question is the hyperbolic Laplacian $\Delta_0$ on $\bbB^{n+1}$, lifted to $Q$
which is the blow-up of $\bbB^{n+1}$ at a point on its boundary.  On each fiber $\sigma$
will be constant, so we can just work with the usual hyperbolic spectral parameter $\zeta$
Let $\rho,\rho'$ be defining functions for the two faces of $Q$, where $\rho=0$ corresponds to the
remnant of the original boundary (these would be the restrictions to the fiber of the coordinates introduced
above).  The space that we are concerned with is
$\calA_{\zeta,\zeta-l}(Q)$ (with logarithmic terms included even though $\zeta$ is constant).
As in (\ref{truncdef}), an index $\zeta|k$ denotes a truncated expansion with leading term
$\rho^{\zeta+k}(\log \rho)^k$.

\begin{proposition}\label{fibermodel}
For $k\in\bbN$, $l\in \bbN_0$, we can extend the model resolvent 
$R_0(\zeta) = [\Delta_0 + \zeta(n-\zeta)]^{-1}$ to a meromorphic map
$$
R_0(\zeta): \calA_{\zeta|k,\zeta-l}(Q) \to  \calA_{\zeta,\zeta-l}(Q),
$$
with poles at $\zeta \in \frac12(n-k-\bbN_0)\cup \frac12(-l-\bbN_0)$ and also
$-\bbN_0$ for $n$ odd.
\end{proposition}
\begin{proof}
It is most convenient to argue in the model for $Q$ given by the upper half-space blown up at the origin.
Sticking to our convention that $x$ is the boundary defining function, we use coordinates $(x,y) \in
\bbR_+\times
\bbR^n$.  The model Laplacian is
$$
\Delta_0 = -(x\del_x)^2 + nx\del_x - (x\del_y)^2.
$$
With radial coordinate $r=\sqrt{x^2+y^2}$,
the defining functions for the faces of $Q$ are $\rho = x/r$ and $\rho'=r$.  According to
Proposition 6.19 of \cite{MM}, the model resolvent extends meromorphically to a map
\begin{equation}\label{619}
R_0(\zeta) = [\Delta_0 - \zeta(n-\zeta)]^{-1}: \rho^{\zeta+k} {\rho'}^{\zeta-l} \cinf(Q) 
\to \rho^\zeta {\rho'}^{\zeta-l} \cinf(Q),
\end{equation}
with poles as indicated above.   

Let $f\in \calA_{\zeta|k,\zeta-l}(Q)$.  We'll solve the problem
\begin{equation}
(\Delta_0 - \zeta(n-\zeta))u = f,
\end{equation}
in stages.  First a simple computation shows that for $\psi \in \calA_{0,\zeta-l}(Q)$,
\begin{equation*}\begin{split}
&(\Delta_0 - \zeta(n-\zeta)) \rho^{\zeta+k} (\log\rho)^m \psi\\
&\qquad= k(n-2\sigma-k) \rho^{\zeta+k} (\log\rho)^m \psi\\
&\qquad\quad+ (m-1)\rho^{\zeta+k} (\log\rho)^{m-1} \psi_1 + (m-1)(m-2)\rho^{\zeta+k} 
(\log\rho)^{m-2} \psi_2 + v,
\end{split}\end{equation*}
where $\psi_1,\psi_2\in \calA_{0,\zeta-l}(Q)$ and $v \in \calA_{\zeta|(k+1),\zeta-l}(Q)$.  We can use this
to solve away the asymptotic expansion in $\rho$, exactly as in Lemma \ref{solveaway}, with poles at 
$\zeta\in\frac12(n-k-\bbN_0)$.  The result is $u_0 \in \calA_{\zeta|k,\zeta-l}(Q)$ such that
$$
(\Delta_0 - \zeta(n-\zeta))u_0 - f = f_1 \in \calA_{\infty,\zeta-l}(Q).
$$

Now suppose we want to solve
\begin{equation}\label{logprob}
(\Delta_0 - \zeta(n-\zeta)) w = (\log\rho') \psi,
\end{equation}
where $\psi\in \rho^{\infty} {\rho'}^{\zeta-l} \cinf(Q)$.  
By (\ref{619}) we can apply $R_0$ to $\psi$.  Then we have
$$
(\Delta_0 - \zeta(n-\zeta)) (\log\rho') R_0(\zeta) \psi = (\log\rho')\psi + [\Delta_0,\log\rho'] R_0(\zeta)
\psi.
$$
A straightforward computation shows that 
$$
[\Delta_0,\log\rho'] \in \rho^2\calV_b(Q),
$$ 
and $\rho^\infty {\rho'}^{\zeta-l} \cinf(Q)$ is invariant under $\calV_b(Q)$.  
Therefore
$$
[\Delta_0,\log\rho'] R_0(\zeta) \psi \in \rho^{\infty} {\rho'}^{\zeta-l} \cinf(Q).
$$
The model resolvent can be applied to this expression by (\ref{619}). The solution to
(\ref{logprob}) is then
$$
w = (\log\rho') R_0(\zeta) f - R_0(\zeta) [\Delta_0,\log\rho'] R_0(\zeta) f,
$$
with poles as indicated.  By induction we can extend this trick to higher powers of $\log\rho'$.

By applying this argument to successive terms in the asymptotic expansion of $f_1$ in $\rho'$
and asymptotically summing the resulting terms, we can find $u_1 \in \calA_{\zeta,\zeta-l}(Q)$
such that
$$
(\Delta_0 - \zeta(n-\zeta))u_1 - f_1 = f_2 \in \dot\cinf(Q).
$$
The solution to the original problem is now given by $u = u_0 - u_1 + R_0(\zeta)f_2$.
\end{proof}

Using Lemma \ref{solveaway} and Proposition \ref{fibermodel}, the parametrix
construction of \cite{MM} can be applied essentially verbatim to the case of variable $\alpha$. 
Note that although a global result is stated in \cite{MM}, the construction is entirely local.
\begin{proposition}\label{param}
For $\zeta\in \bbC/\Gamma$ there is a parametrix $M_\zeta$, analytic in $\zeta$, such that
$$
\Dz M_\zeta = I - F_\zeta
$$
where $M_\zeta \in \psop^{-2}(X) + \psop_{\sigma_l,\sigma_r}(X)$ and the error term $F_\zeta \in
\Psi_{\infty,\sigma_r}(X)$.  

If $\zeta\in \Gamma$ we can  construct a local parametrix
with the same structure, defined in a neighborhood $G\times G \subset \xx$ such that $\sigma(\zeta,y)\notin
\frac12 (n-\bbN_0)$ on $G\cap \bX$.
\end{proposition}

\section{Partial meromorphic continuation}\label{merosec}

Consider the parametrix $M_\zeta$ of Proposition \ref{param} for $\zeta\in
\bbC\bck\Gamma$.    As in the asymptotically hyperbolic case the error $F_\zeta$ is
compact on weighted $L^2$ spaces so by the analytic Fredholm theorem there is a
meromorphic inverse $(I-F_\zeta)^{-1}$.  The resolvent is then given by
$M_\zeta(I-F_\zeta)^{-1}$.  Let $I+D_\zeta = (I-F_\zeta)^{-1}$.  Then the relations
$$
D_\zeta = F_\zeta + D_\zeta F_\zeta = F_\zeta + F_\zeta D_\zeta
$$
can be used to show that $D_\zeta \in \Psi_{\infty,\sigma_r}(X)$ also.  We then claim that
$M_\zeta D_\zeta \in \Psi_{\sigma_l,\sigma_r}(X)$.  Since this amounts to a special case of Lemma
\ref{residlem}, we will not give a separate proof here.  The result is:

\begin{theorem}\label{merores}
As an operator $\dot\cinf(X)\to \cinf(X^{\circ})$, the resolvent $R_\zeta = [\Delta -
\alpha_0^2\zeta(n-\zeta)]^{-1}$ has a meromorphic continuation to $\zeta\in\bbC\bck\Gamma$.  Moreover,
the resolvent has the structure
$$
R_\zeta \in \psop^{-2}(X) + \psop_{\sigma_l,\sigma_r}(X) + \Psi_{\sigma_l,\sigma_r}(X).
$$
\end{theorem}

As discussed in the introduction, in the asymptotically hyperbolic case the existence of a Poisson kernel
and scattering matrix follows directly from this theorem \cite{JS}.  We'll prove the corresponding
implications of Theorem \ref{merores} in this section.  The proofs are similar to those of \cite{JS},
with some modifications necessitated by the variable orders.

As in the asymptotically hyperbolic case, the Poisson kernel will be obtained by
restriction of the resolvent kernel.  We continue to use the Riemannian density of $g$ to identify
operators with integral kernels.  With this convention the Poisson kernel will be given up to a constant
by
\begin{equation}
E_\zeta = {x'}^{-\sigma_r} R_\zeta|_{x'=0}.
\end{equation}
$E_\zeta$ is most naturally described by its lift to 
$$
\sxbx = [\xbx; S],
$$
where $S = \Delta(\bX\times\bX)$ as before.  This $\sxbx$ is naturally diffeomorphic to
the right face of $\sxx$ and has two boundary hypersurfaces, left and front.

Let $R_\zeta = A_\zeta + B_\zeta + C_\zeta$, decomposed as in Theorem \ref{merores}.
Since $A_\zeta$ is supported in a neighborhood of the lifted diagonal, we have
$$
{x'}^{-\sigma_r} A_\zeta|_{x'=0} = 0.
$$
Thus $E_\zeta = H_\zeta+Q_\zeta$, where
$$
H_\zeta = {x'}^{-\sigma_r} B_\zeta|_{x'=0} \in \calA_{\sigma_l, -\sigma_l}(\sxbx),
$$
and 
$$
Q_\zeta = {x'}^{-\sigma_r} C_\zeta|_{x'=0} \in \calA_\sigma(\xbx).
$$
The characterization of $H_\zeta$ and $Q_\zeta$ as polyhomogeneous with variable order follows immediately
from Proposition \ref{bexp}.

Let $h_0 = h(0,y,dy)$ be the metric on $\bX$ induced by $\bar g$.
For $f\in \bX$, define 
\begin{equation}\label{ezfdef}
E_\zeta f(x,y) = \int_\bX E_\zeta(x,y,y') f(y')\>dh_0(y').
\end{equation}
We will show that $E_\zeta f$ has an asymptotic expansion for $x=0$ and compute the leading terms.

\begin{lemma}\label{dblseries}
Let $u\in\cinf(X^{\circ})$ and $\alpha,\beta\in\cinf(X)$ with $\alpha-\beta \notin \bbZ$.  
Suppose that for every $m>0$ there exists $q\in\bbN$ such that
$$
\prod_{k=0}^q\Bigl[x^{\alpha}(x\del_x-k)^{k+1}x^{-\alpha}\Bigr]
\prod_{k=0}^q\Bigl[x^{\beta}(x\del_x-k)^{k+1}x^{-\beta}\Bigr]u
\in \calA^{m-}(X).
$$
Then
$$
u\in\calA_\alpha(X) + \calA_\beta(X).
$$
\end{lemma}
\begin{proof}
Let
\begin{equation}\label{updef}
u_p = \prod_{k=p+1}^q\Bigl[x^{\alpha}(x\del_x-k)^{k+1}x^{-\alpha}\Bigr]
\prod_{k=0}^q\Bigl[x^{\beta}(x\del_x-k)^{k+1}x^{-\beta}\Bigr]u,
\end{equation}
so that
$$
\prod_{k=0}^{p} (x\del_x-k)^{k+1} (x^{-\alpha}u_p) \in x^{-\alpha}\calA^{m-}(X)
$$
Assuming that $p<m-\re\alpha$, we can apply the argument from Lemma \ref{bexp} to get
$a_{j,l}\in\cinf(X)$ such that
\begin{equation}\label{upexp}
x^{-\alpha} u_p - \sum_{j=0}^{p-1}\sum_{l=0}^j x^j(\log x)^l a_{j,l} \in \calA^{p-}(X).
\end{equation}

Now If we define
\begin{equation}\label{wdef}
w = \prod_{k=0}^q\Bigl[x^{\beta}(x\del_x-k)^{k+1}x^{-\beta}\Bigr]u,
\end{equation}
then from (\ref{updef}) the relation to $u_p$ is
$$
x^{-\alpha} u_p = \prod_{k=p+1}^q (x\del_x-k)^{k+1} (x^{-\alpha} w).
$$
By solving a linear system for the coefficients $b_{j,l}$ in terms of the $a_{j,l}$
we can rewrite (\ref{upexp}) as
\begin{equation}\label{kpq}
\prod_{k=p+1}^q (x\del_x-k)^{k+1} \Bigl(x^{-\alpha} w
- \sum_{j=0}^{p-1}\sum_{l=0}^j x^{j}(\log x)^l b_{j,l}\Bigr) \in \calA^{p-}(X)
\end{equation}

Suppose $F\in \cinf(X^\circ)$ satisfies
\begin{equation}\label{xxkF}
(x\del_x-k) F \in \calA^{p-}(X),
\end{equation}
for $k>p$.  This implies
$$
\del_x (x^{-k} F) = h\in x^{-k-1} \calA^{p-}(X).
$$
Integrating from $x$ to $1$ gives
$$
x^{-k} F = F(1) - \int_x^1 h\>dx.
$$
We note that $(\calV_b)^l \int_x^1 h\>dx = O(x^{p-k-\epsilon})$ for $\epsilon>0$, so we conclude that
(\ref{xxkF}) implies $F\in \calA^{p-}(X)$.
Applying this repeatedly to (\ref{kpq}) gives
$$
w - \sum_{j=0}^{p-1}\sum_{l=0}^j x^{\alpha+j}(\log x)^l b_{j,l} \in x^{\alpha}\calA^{p-}(X).
$$
 
Substituting back with the definition of $w$ from (\ref{wdef}) and once again solving a linear system for
new coefficients, we get
\begin{equation}\label{fstpass}
\prod_{k=0}^q\Bigl[x^{\beta}(x\del_x-k)^{k+1}x^{-\beta}\Bigr]
\Bigl(u - \sum_{j=0}^{p-1}\sum_{l=0}^j x^{\alpha+j} (\log x)^l c_{j,l}\Bigr) \in x^{\alpha}\calA^{p-}(X)
\end{equation}
The matrix relating the $c_{j,l}$'s to the $b_{j,l}$'s is lower triangular, with diagonal
entries of the form $(\alpha+j-\beta-k)$.  Hence the requirement that $\alpha-\beta\notin \bbZ$
ensures the system is non-singular.

Now we let $v$ be the function in parentheses in (\ref{fstpass}) and
simply repeat the argument given above.  Assuming
$s<p+\re\alpha-\re\beta$,  we obtain $d_{j,l}\in \cinf(X)$ such that
$$
u - \sum_{j=0}^{p-1}\sum_{l=0}^j x^{\alpha+j} (\log x)^l c_{j,l}
- \sum_{j=0}^{s-1}\sum_{l=0}^j x^{\beta+j} (\log x)^l d_{j,l} \in x^{\beta}\calA^{s-}(X).
$$
This construction is possible for $s + \re\beta < p + \re\alpha < m$, and $m$ can be arbitrarily large,
so we have full asymptotic expansions.
\end{proof}

\begin{proposition}\label{ezfprop}
For $f\in \dot\cinf(W_\lambda)$ and $E_\zeta f$ defined by (\ref{ezfdef}) we have
$$
E_\zeta f \in \calA_\sigma(X) + \calA_{n-\sigma}(X).
$$
\end{proposition}
\begin{proof}
That $Q_\zeta f\in\calA_{\sigma}(X)$ follows immediately from $Q_\zeta\in\calA_\sigma(\xbx)$,
so we concentrate on $H_\zeta$.  The asymptotic properties of $E_\zeta f$ depend only on 
the behavior of $E_\zeta(x,y,y')$ near $x=0$ and/or $y=y'$.  So we can specialize to a particular
coordinate neighborhood with coordinates $(x,y,z=y-y')$ and assume that all functions are compactly
supported within this neighborhood.  We can rewrite (\ref{ezfdef}) as
$$
E_\zeta f(x,y) = \int_{\bbR^n} w(x,y,y-z) \>dz,
$$
for $w\in\calA_{\sigma_l, -\sigma_l}(\sxbx)$.  On $\sxbx$ we'll use coordinates $r = \sqrt{x^2+z^2}$,
$\rho = x/r$, $\omega = z/r$, and $y$.  For convenience, we extend $\sigma$ into the interior so as to
be independent of $x$ within the neighborhood of interest.  So $\sigma_l = \sigma(y)$ does not
depend on $r$, $\rho$, or $\omega$.

Noting that
$$
\int z\del_z w\>dz = -n \int w \>dz
$$
by integration by parts, we can move derivatives under the integral to get
\begin{equation}\label{dui}
(x\del_x-\sigma)(x\del_x - n+\sigma)E_\zeta f
= \int (x\del_x-\sigma_l)(x\del_x +z\del_z+\sigma_l) w\>dz
\end{equation}
The lift of the vector fields appearing here to $\sxbx$ is
\begin{equation*}\begin{split}
x\del_x &\longrightarrow \rho^2 r\del_r + (1-\rho^2) \rho\del_\rho + \rho^2 \omega\del_\omega\\
x\del_x +z\del_z &\longrightarrow r\del_r
\end{split}\end{equation*}
Writing $w = \rho^{\sigma_l} r^{-\sigma_l} F$ for $F\in \calA_{0,0}(\sxbx)$, the integrand
on the right in (\ref{dui}) becomes
$$
(x\del_x-\sigma_l)(r\del_r + \sigma_l) w =  \rho^{\sigma_l} r^{-\sigma_l} 
\bigl[\rho\del_\rho + \rho^2(-2\sigma_l + r\del_r - \rho\del_\rho +
\omega\del_\omega)\bigr](r\del_r) F 
$$
Now by the definition of the truncated spaces in \S\ref{basymsec}),
\begin{equation*}\begin{split}
(r\del_r)F &\in  \calA_{0,0|1}(\sxbx);\\
(\rho\del_\rho)(r\del_r)F &\in \calA_{0|1,0|1}(\sxbx).
\end{split}\end{equation*}
Also $r\del_r - \rho\del_\rho + \omega\del_\omega \in \calV_b(\sxbx)$ and 
$$
\rho^2 \calA_{0,0|1}(\sxbx)\subset \calA_{0|2,0|1}(\sxbx),
$$   
so we conclude that
$$
(x\del_x-\sigma_l)(r\del_r + \sigma_l) w \in \rho^{\sigma_l} r^{-\sigma_l}\calA_{0|1,0|1}
= \calA_{\sigma_l|1,-\sigma_l|1}
$$

Applying this argument inductively gives
\begin{equation}
\prod_{k=0}^{q-1}\Bigl[(x\del_x-\sigma-k)(x\del_x-n+\sigma-k)\Bigr]^{k+1} E_\zeta f
= \int G\>dz,
\end{equation}
where $G \in \calA_{\sigma_l|q,-\sigma_l|q}$.  For $q > m+\re\sigma$ we have
$\calA_{\sigma_l|q,-\sigma_l|q} \subset \calA^{(m,m)-}(\sxbx)$.  The lift of $\calV_b(X)$ to
$\sxbx$ through the left lies in $r^{-1}\calV_b(\sxbx)$, so by passing derivatives under the integral and
using the estimates on $G$ we have
$$
\int G\>dz \in \calA^{(m-1)-}(X).
$$
Hence for $q > m+\re\sigma$ we obtain
$$
\prod_{k=0}^{q-1}\Bigl[(x\del_x-\sigma-k)(x\del_x-n+\sigma-k)\Bigr]^{k+1} E_\zeta f
\in \calA^{(m-1)-}(X).
$$
This has been derived under the assumption that $\sigma$ is independent of $x$ in the local neighborhood,
so the result then follows from Lemma \ref{dblseries}.
\end{proof}

Proposition \ref{ezfprop} shows that $E_\zeta f$ has two separate asymptotic expansions at
$x=0$.  Let $x^{\sigma}a_0$ and $x^{n-\sigma}b_0$ be the leading terms for each. 
The coefficients $a_0$ and $b_0$ will be holomorphic functions of $\zeta$ for $\zeta\in\bbC\bck\Gamma$,
which is very useful for computing them.
If $\zeta$ is such that $2\re\sigma - n >0$, then
$$
b_{0}(\zeta,y) = \lim_{x\to 0} x^{\sigma-n} \int_{\bbR^n} E_\zeta(x,y,y') \phi(y')\>dy',
$$
where $\phi = f \sqrt{h_0}$.
Introducing the coordinate $w = (y-y')/x$, we have
$$
b_{0}(\zeta,y) = \lim_{x\to 0} \int_{\bbR^n} x^{\sigma} E_\zeta(x,y,y-xw) \phi(y-xw) \>dw.
$$
Note that $x^{\sigma}E_\zeta \in \calA_{2\sigma_l,0}$.  Write this as $\rho^\sigma F(r,\rho,\omega)$, where
$F\in \calA_{0,0}$.  Noting that $r = x\sqrt{1+|w|^2}$, $\rho = 1/\sqrt{1+|w|^2}$ and 
$\omega = \rho w$, we have
$$
b_{0}(\zeta,y) = \phi(y) \int_{\bbR^n}  F(0,\rho,\rho w) \frac{dw}{(1+|w|^2)^{\sigma}}. 
$$
Since $F$ is bounded and $2\re\sigma>n$, the integral is convergent.

As pointed out in \cite{JS}, $F(0,\rho,\rho w)$ can be computed from the restriction of 
$R_\zeta$ to the front face, i.e. the normal operator of $R_\zeta$.  Since this is just the model
resolvent, we see that
\begin{equation}\label{b0}
b_{0} = B(\sigma) f \sqrt{h_0},
\end{equation}
where $B$ is a universal meromorphic function depending on $n$ but not on $g$.  This formula is
therefore valid for all $\zeta\in \bbC\bck\Gamma$.  Combining Proposition \ref{ezfprop} with (\ref{b0})
gives us the following result: for $f\in\cinf(\bX)$ the function 
$$
u = \frac{1}{B(\sigma) \sqrt{h_0}} E_\zeta f
$$
satisfies $\Dz u = 0$ and
$$
u \sim x^{n-\sigma} f + x^\sigma f',
$$
with $f'\in\cinf(\bX)$.
We have thus proven:
\begin{proposition}\label{poissk}
With respect to the metric density of $h_0$ the Poisson kernel is
$$
\frac{1}{B(\sigma) \sqrt{h_0}} E_\zeta.
$$
\end{proposition}

To understand how $f'$ relates to $f$ we must compute $a_0$. 
If $2\re\sigma -n<0$ then 
$$
a_{0}(\zeta,y) = \lim_{x\to 0} x^{-\sigma} \int_{\bbR^n} E_\zeta(x,y,y-z) \phi(y-z)\>dz.
$$
Using the same definition of $F$ as above, 
$$
x^{-\sigma} E_\zeta = r^{-2\sigma} F(r,\rho,\omega).
$$
The limit can be taken directly
$$
a_{0}(\zeta,y) = \int_{\bbR^n} |z|^{-2\sigma} F(|z|,0,z/|z|) \phi(y-z)\>dz.
$$
Hence the scattering matrix is given by 
$$
S_\zeta:f \mapsto \frac{1}{B(\sigma) \sqrt{h_0}}
\int_{\bbR^n} \frac{F(|y-y'|,0,\frac{y-y'}{|y-y'|})}{|y-y'|^{-2\sigma}} f(y')\>dh_0(y')
$$
In computing the symbol of this pseudodifferential operator, the two factors of $\sqrt{h_0}$
cancel out (they would not have appeared if we worked with half-densities).
\begin{proposition}\label{smatrix}
For $\zeta\in\bbC\bck\Gamma$ the scattering matrix $S_\zeta$ is a meromorphic pseudodifferential
operator on $\bX$, with principal symbol 
$$
2^{n-2\sigma} \frac{\Gamma(\frac{n}2-\sigma)}{\Gamma(\sigma-\frac{n}2)} |\xi|^{2\sigma-n}_{h_0}.
$$
\end{proposition}

\section{Limiting absorption}\label{lasec}

Let $R_\zeta = \Dz^{-1}$, as defined by the partial meromorphic continuation. 
Fix $\zeta$ with $\re\zeta=n/2$ and $\im\zeta\ne 0$, and choose a sequence
$\zeta_j \to\zeta$ with $\re\lambda_j>n/2$ and $\im\zeta_j = \im\zeta$.  
For fixed $f \in \dot\cinf(X)$, the point of this section to show
that the sequence $R_{\zeta_j}f$ converges to the unique solution $u$ of $\Dz u = f$
satisfying a certain `radiation' condition.  By this technique we are able to define $R_\zeta$
for $\zeta$ on the boundary of the region of analytic continuation, corresponding to the irregular
continuous spectrum $\lambda\in (\frac{\alpha_0^2n^2}4, \frac{\alpha_1^2n^2}4]$. 
We have to avoid $\lambda=\frac{\alpha_0^2n^2}4$ because of the possible
embedded eigenvalue there.  However, if we had $\alpha=\alpha_0$ on some open set in $\bX$ 
then an embedded eigenvalue would be ruled out and the point $\zeta=\frac{n}2$ could be
included as a possible limit point of $\{\zeta_j\}$.

The limiting absorption property relies on certain uniform estimates on $R_{\zeta}f$ for
$\re\zeta>\frac{n}2$.  
In the following, for complex valued 1-forms $\theta_1,\theta_2$ we'll denote the metric
inner product by
$$
\pair{\theta_1,\theta_2} = g^{-1}(\ol{\theta_1},\theta_2).
$$
For either functions or forms, $(\cdot,\cdot)$ denotes the $L^2$ inner product, and $\norm\cdot$ the 
$L^2$ norm.

To control the boundary terms in the propositions below, we first state a non-uniform result which is a
simple consequence of Theorem \ref{merores} and Lemma \ref{maplemma}.
\begin{lemma}\label{rfest}
For $\re\zeta > n/2$, $f\in\dot\cinf(X)$, there exists $\epsilon>0$ (depending on $\zeta$) such that
$$
R_\zeta f = O(x^{n/2+\epsilon}).
$$
\end{lemma}

The first uniform estimate is relatively simple.
\begin{proposition}\label{normdu}
Let $u = R_\zeta f$ for $f\in\dot\cinf(X)$ and $\re\zeta>n/2$.  For any $\delta>0$
$$
\norm{x^\delta du}^2 \le C\Bigl(|\lambda|\>\norm{x^\delta u}^2 + \norm{f}^2\Bigr),
$$
where $\lambda = \alpha_0^2\zeta(n-\zeta)$ and $C$ is independent of $\zeta$.
\end{proposition}
\begin{proof}
Integrating by parts gives
$$
\int_{x\ge s} x^{2\delta} \ol{u} \Delta u \>dg = \int_{x\ge s} \pair{d(x^{2\delta}
\ol{u}),  du} \>dg + \int_{x=s} x^{2\delta} \ol{u} (-\alpha x\del_x u) \frac{dh}{x^n}
$$
In taking the limit $s\to\infty$, the boundary term disappears by Lemma \ref{rfest},
yielding
\begin{equation}
\int x^{2\delta} \ol{u} \Delta u \>dg = \int \pair{d(x^{2\delta} \ol{u}),  du} \>dg.
\end{equation}
Using $(\Dl) u = f$, we rewrite this as
\begin{equation*}
\begin{split}
\norm{x^\delta du}^2 &= \int x^{2\delta} \ol{u} (\lambda u + f) \>dg - 2\delta \int
x^{2\delta-1} \ol{u} \pair{u\> dx, du}\> dg\\ 
&= \lambda \norm{u}^2 + (x^{2\delta}u,f) - 2\delta (x^{\delta-1} u\>dx, 
x^{\delta} du)
\end{split}
\end{equation*}
Thus we can estimate
$$
\norm{x^\delta du}^2 \le C\Bigl(\norm{x^\delta u}^2 + \norm{x^{2\delta} u}\>\norm{f}
+ \norm{x^\delta u}\> \norm{x^\delta du}\Bigr).
$$
We take care of the $\norm{x^\delta du}$ on the right by an estimate of the form
$$
\norm{x^\delta u}\> \norm{x^\delta du} \le \epsilon \norm{x^\delta du}^2 +
\frac1{4\epsilon} \norm{x^\delta u},
$$
and the result follows.
\end{proof}

The second estimate establishes the radiation condition.  It is also essentially just integration by parts,
but a more complicated computation.  The strategy is taken from a similar argument by Tayoshi 
\cite{Ta}, who established limiting absorption for a class of asymptotically flat metrics.

Fix $\zeta_0\ne \frac{n}2$ with $\re\zeta_0=n/2$ and let $I_{\zeta_0} = \{\zeta_0+t: \; 0<t<1\}$.    
For $\zeta\in I_{\zeta_0}$, the indicial root $\sigma$ is smooth, but not uniformly so as $\zeta\to \zeta_0$.
We thus need to be a little careful about the extension of $\sigma$ into the interior.
Away from the set $\{\alpha^2 = 4\lambda/n^2\}$ where $\sigma(\zeta_0)$ is singular any smooth extension will
do.  Near this set, however, we need to impose some uniformity.  In such neighborhoods we will take
$$
\tsigma(\zeta,x,y) = \frac{n}2 + \sqrt{\frac{n^2}4 - \frac{\alpha_0^2}{\alpha^2} \zeta(n-\zeta) + ix}.
$$
Such an extension has the following properties:
\begin{enumerate}
\item $\tsigma \in C^0(X)$ is smooth in the interior and $\tsigma|_\bX = \sigma$.
\item The difference $\alpha^2 \tsigma(n-\tsigma) - \alpha_0^2 \zeta(n-\zeta)$ is $O(x)$,
uniformly for $\zeta\in I_{\zeta_0}$.
\item $|d\tsigma| = O(x^{1/2})$, uniformly for $\zeta\in I_{\zeta_0}$.
\end{enumerate}

\begin{proposition}\label{radiation} 
Fix $f\in\dot\cinf(X)$ and extend $\sigma$ into the interior
as above.  If $u = R_\zeta f$ for $\zeta\in I_{\zeta_0}$ then we can estimate   
$$
\Bigl\Vert x^{-\epsilon} \bigl(du - \tsigma(\zeta) u \>\tfrac{dx}x\bigr)
\Bigr\Vert^2 < C (\norm{x^\delta u}^2 + \norm{f}^2),
$$
where $\delta, \epsilon>0$, $\delta+\epsilon<1/2$, and $C$ is independent of $\zeta$. 
\end{proposition}

\begin{proof}
Let us single out the operator appearing in the radiation condition by defining
$$
\ds u = du - \tsigma u \tfrac{dx}{x}.
$$
We first note that the only issue is to control $\ds u$ at the boundary.  Because for $t>0$,
$$
\int_{x\ge t} |\ds u|^2 \>dg \le C \int x^{2\delta} |\ds u|^2 \>dg,
$$
for any $\delta>0$, and by Proposition \ref{normdu}, we have
$$
\int x^{2\delta} |\ds u|^2 \>dg \le C (\norm{x^\delta u}^2 + \norm{f}^2)
$$
So we will work in a product neighborhood of $\bX$ in which the metric has the form (\ref{norform}).
Since $|\frac{dx}{\alpha x}| =1$, the radial component of $\ds u$ can be singled out as
$$
\Ds u = \pair{\tfrac{dx}{\alpha x}, \ds u} = \alpha(x\del_x u - \tsigma u).
$$
We will also let $\gamma = x\del_x \log \sqrt{h}$ so that
$$
\dv(\del_x) = \frac{n+1-\gamma}{x}.
$$

Let $s>0$ be small enough that $[0,s]\times\bX$ lies within the product neighborhood.
For $\epsilon>0$ we choose some $\psi\in \cinf(0,\infty)$ such that
$$
\psi(x) = \begin{cases}x^{-2\epsilon}& \text{if }x\le s/2 \\ 0&\text{if } x>s\end{cases}
$$

We begin by computing a divergence:
\begin{equation}\label{divdd}
\dv [\alpha^{-1}\psi \ol{\Ds u} \>g^{-1}(d_\tsigma u)] =  \alpha^{-1}\psi \ol{\Ds u}
\> \dv g^{-1}(\ds u) + \pair{d(\alpha^{-1}\psi\Ds u), \ds u} 
\end{equation}
In the first term on the right-hand side we have
\begin{equation*}\begin{split}
\dv g^{-1}(\ds u) &= -\Delta u + \dv (\alpha^2 \tsigma u x \del_x)\\
&= -\alpha_0^2 \zeta(n-\zeta) u - f + \alpha^2\tsigma x\del_x u + \alpha^2\tsigma (n-\gamma) u + \alpha^2 (x\del_x
\tsigma) u\\
&= -f + \alpha\tsigma \Ds u + Ru,
\end{split}\end{equation*}
where
$$
R = \alpha_0^2\zeta(n-\zeta) - \alpha^2\tsigma(n-\tsigma) + \alpha^2(\tsigma\gamma - x\del_x \tsigma).
$$
Inserting back into (\ref{divdd}), we have
\begin{equation}\label{dveq}\begin{split}
&\dv [\alpha^{-1}\psi \ol{\Ds u} \>g^{-1}(d_\tsigma u)]\\ 
&\quad= \pair{d(\alpha^{-1}\psi\Ds u),\ds u}
+ \tsigma\psi|\Ds u|^2 + \alpha^{-1}\psi\ol{\Ds u} (Ru-f)\\
&\quad= \psi \pair{\ds(\alpha^{-1}\Ds u),\ds u} + x \psi' |\Ds u|^2 
+ \alpha^{-1}\psi\ol{\Ds u} (Ru-f)
\end{split}\end{equation}
Now integrating the divergence on the left-hand side of (\ref{dveq}) gives
$$
\int \dv [\alpha^{-1}\psi \ol{\Ds u} \>g^{-1}(d_\tsigma u)]  dg 
= \lim_{t\to 0^+} \int_{x=t} \alpha^{-1} \psi |\Ds u|^2 \>\frac{dh}{x^n},
$$
and by Lemma \ref{rfest} the boundary term is zero.  Hence the integral of (\ref{dveq}) gives 
\begin{equation}\label{divint}
0=\int \Bigl[\psi \pair{\ds(\alpha^{-1}\Ds u),\ds u} + x \psi' |\Ds u|^2 
+ \alpha^{-1}\psi\ol{\Ds u} (Ru-f) \Bigr] dg
\end{equation}

One more integration by parts is needed, to get rid of the second derivative of $u$ appearing in the 
first term of the integrand in (\ref{divint}).  To this end, note that
\begin{equation}\label{dsd}
\begin{split}
\ds(\alpha^{-1}\Ds u) &= \ds(x\del_x u - \tsigma u) \\
&= \calL_{x\del_x} (\ds u) - u \>d\tsigma - \tsigma \ds u.
\end{split}
\end{equation}
On substituting this into (\ref{divint}), the second derivative of $u$ would be contained in the term
$$
\int \psi \pair{\calL_{x\del_x} \ds u, \ds u}
$$
This can be integrated by parts to give 
$$
- \frac12 \int \Bigl[(x \psi' + (\gamma- n)\psi) |\ds u|^2 + \psi  
(\calL_{x\del_x} g^{-1})(\ol{\ds u},\ds u) \Bigr] dg,
$$
with a boundary term of 
$$
\lim_{t\to 0^+} \frac12 \int_{x=t} \alpha^{-1} \psi |\ds u|^2 \>\frac{dh}{x^n} = 0
$$
(again by Lemma \ref{rfest}).
Combining this integration by parts with (\ref{dsd}) and substituting back into (\ref{divint}) we now have
\begin{equation}\label{intC}
\begin{split}
0 = \int & \Bigl[ 
- \tfrac12 (x \psi' + (\gamma - n)\psi) |\ds u|^2 - \tfrac12 \psi (\calL_{x\del_x} g^{-1})(\ol{\ds
u},\ds u) \\  
&- \psi \pair{u \>d\tsigma, \ds u} - \psi\tsigma |\ds u|^2 + x \psi' |\Ds u|^2  +  \alpha^{-1}\ol{\Ds u}\>
(R u + f)
 \Bigr] dg
\end{split}
\end{equation}

The Lie derivative of $g^{-1}$ may be expressed in terms of $h^{-1}$: 
$$
- \tfrac12 \calL_{x\del_x} g^{-1} = x^2h^{-1} + \tfrac12 x^3\del_x h^{-1}
$$
Since $g^{-1} = \alpha^2 x^2 \del_x\otimes \del_x + x^2h^{-1}$, for
some $k<1$ fixed we will have
\begin{equation}\label{lieg}
- \tfrac12 \calL_{x\del_x} g^{-1} - k [g^{-1} - \alpha^2 x^2 \del_x\otimes \del_x] \ge 0,
\end{equation}
as a tensor, for $x$ sufficiently small.  Recalling that $\psi$ was to be supported in $[0,s]$,
we require that $s$ be sufficiently small so that (\ref{lieg})
holds for $x<s$.  Applying (\ref{lieg}) to the relevant term in (\ref{intC})  gives
$$
- \tfrac12 \psi (\calL_{x\del_x} g^{-1})(\ol{\ds u},\ds u) \ge  k \psi (|\ds u|^2 - |\Ds u|^2)  
$$

Using this in (\ref{intC}) along with the fact that
$\re\tsigma>\frac{n}2$, we obtain an inequality
\begin{equation}\label{intD}
\begin{split}
&\int  \Bigl[
(\tfrac12 x \psi' + \tfrac12 \psi \gamma + k\psi) |\ds u|^2 - (x \psi'+k\psi) |\Ds u|^2\Bigr] dg\\
&\quad\le \int \re\Bigl[
- \psi \pair{u \>d\tsigma, \ds u} +  \psi \ol{\Ds u}\> (R u + \alpha^{-1} f) \Bigr] dg
\end{split}
\end{equation}

Recalling that $\psi=x^{-2\epsilon}$ for $x<s/2$, we split the integral on the left-hand side of
(\ref{intD}) to yield 
\begin{equation}\label{ineqA}\begin{split}
&\int_0^{s/2}  x^{-2\epsilon} \Bigl[
(-\epsilon +  \tfrac12 \gamma + k) |\ds u|^2  - (k-2\epsilon)|\Ds u|^2 \Bigr] dg \\
&\le \int x^{-2\epsilon} \re\Bigl[
-  \pair{u \>d\tsigma, \ds u} + \ol{\Ds u}\> (R u + \alpha^{-1} f) \Bigr] dg
 + C \int_{s/2}^\infty |\ds u|^2 dg
\end{split}\end{equation}
We can pick $\epsilon <k/2$, and since $\gamma = O(x)$, 
by choosing $s$ sufficiently small we can assume $\epsilon + \gamma/2 > c > 0$.
Noting also that $|\Ds u|\le |\ds u|$, the left-hand side can be bounded below by
$$
c\int_0^{s/2}  x^{-2\epsilon} |\ds u|^2 \> dg. 
$$
Applying this bound together with the estimates
$$
|d\tsigma| = O(x^{1/2}), \qquad R = O(x^{1/2}),
$$
we derive from (\ref{ineqA}) the new inequality
\begin{equation}\label{intE}
\begin{split}
\int_0^{s/2} x^{-2\epsilon} |\ds u|^2 dg \le  C &\int  \Bigl[
x^{1/2-2\epsilon} |u| \>|\ds u| + x^{1/2-2\epsilon}
|\Ds u|\> |u|\\ &\qquad + x^{-2\epsilon}|\Ds u| \> |f|\Bigr] dg
+ C \int_{s/2}^\infty |\ds u|^2 dg.
\end{split}
\end{equation}
By Proposition \ref{normdu} the $L^2$ norm of $x^{1+\delta}\del_y u$ is bounded by that of
$x^\delta u$ for $\delta>0$.  And 
$\int_{s/2}^\infty |\ds u|^2 \> dg$ is estimated as indicated at the start of the proof.
With these facts we can reduce (\ref{intE}) to
$$
\norm{x^{-\epsilon} \ds u}^2 \le C \Bigl(\norm{x^\delta u}^2 + \norm{x^{1/2-\epsilon} u}
\norm{x^{-\epsilon} \Ds u} + \norm{f} \norm{x^{-\epsilon} \Ds u}\Bigr).
$$
Since $|\Ds u| \le |\ds u|$, if we set $\delta<1/2-\epsilon$ this yields the stated inequality.
\end{proof}

The final step in establishing limiting absorption is a uniqueness result needed to guarantee that the
sequence $R_{\zeta_j}f$ converges.

\begin{proposition}\label{unique}
Fix $\lambda>\frac{n^2\alpha_0^2}4$ and let $G\subset X$ be some neighborhood intersecting $\bX$ in 
$W_\lambda$.  Suppose on $G$ there exists a solution $u$ to
$$
[\Dl] u = 0, \qquad u \in x^{-\delta}L^2(G)
$$
for $\delta>0$, satisfying the radiation condition: 
$$
\ds u \in x^\gamma L^2(G; T^*X),
$$
for $\gamma>0$.  Then $u=0$.
\end{proposition}
\begin{proof}
Mazzeo argued in \cite{Ma91a} that the existence of a local parametrix in $G\times G$ as in
Proposition \ref{param} shows that
there can be no solution of $[\Dl] u = 0$ lying in $L^2(G)$. So the point here is
to use the radiation condition to argue that $u$ would have to be $L^2$.

Let $\chi$ be a cutoff in $\cinf_0(G\cap \bX)$.
Choose $\psi\in \cinf(\bbR_+)$ with $\psi = 0$ for $x\le 1/4$, $\psi = 1$ for $x \ge 2$, and
$\psi' \ge 1$ for $1/2 \le x\le 1$.  Setting $\pse(x) = \psi(x/\epsilon)$, we can interpret $\pse \chi$
as a function on $X$ using the product neighborhood on which the metric has the normal form
(\ref{norform}).  Assume $\epsilon$ is small enough that $\pse \chi$ is supported in $G$.

Assuming that $u$ is as stated, we have
\begin{equation}\label{udpse}
(u, [\Delta, \pse \chi]u) = (\ol\lambda - \lambda) (u, \pse \chi u) = 0
\end{equation}
Choose $\epsilon$ small enough that $[\Delta, \pse \chi]u=0$ for $x\ge 1$.
Then (\ref{udpse}) becomes
\begin{equation}\label{commint}
\int_0^1 \ol{u} [\Delta, \pse \chi]u \>dg = 0.
\end{equation}

The commutator is
$$
[\Delta, \pse \chi]u = (\chi \Delta \pse + \pse \Delta \chi) - 2\chi (x\del_x \pse)
(x\del_x u) - 2\pse x^2 h^{ij} (\del_j \chi)(\del_j u) 
$$
Substituting this back into (\ref{commint}) and taking the imaginary part yields
\begin{equation}\label{impart}
2\im \int_0^1 \chi \ol{u} (x\del_x \pse) (x\del_x u) \>dg = - 2\im\int_0^1 \ol{u} \pse x^2
h^{ij} (\del_j \chi)(\del_j u) \>dg
\end{equation}
Let $v = (x\del_x - \sigma)u$, which lies in $x^\gamma L^2(G)$ by assumption.  Then we can write
(\ref{impart}) as
\begin{equation}\label{intpse}
\int_0^1 \chi |u|^2 (\im\sigma) (x\del_x \pse) \>dg = - \im\int_0^1 \ol{u} \Bigl[\pse x^2
h^{ij} (\del_j \chi)(\del_j u) + \chi (x\del_x \pse) v \Bigr]\>dg
\end{equation}
The first term on the right can be bounded:
$$
\Bigl| \int_0^1 \ol{u} \Bigl[\pse x^2 h^{ij} (\del_j \chi)(\del_j u) \>dg \Bigr|
\le C \int_{\epsilon/4}^{2\epsilon} x\>|u|\> |d_y u| \>dg
$$
By assumption $u\in x^{-\delta}L^2$ and $d_y u \in x^{\gamma}L^2$, so $|u| |d_y u|$
is integrable with respect to $dg$.  Thus this first term can be
estimated $O(\epsilon)$.  For the second term on the right-hand side of (\ref{intpse})
we have the bound
$$
 \Bigl| \int_0^1 \chi (x\del_x \pse) \ol{u} v \>dg \Bigr| \le C \int_{\epsilon/4}^{2\epsilon}
\frac{x}\epsilon |u|\>|v|\>dg.
$$
Since $x^{-\gamma+\delta}|u|\>|v|$ is integrable, this term can be
estimated $O(\epsilon^\gamma)$ for small $\gamma>0$.

Since $\im \sigma$ is bounded from below in the support of $\chi$, using
these estimates back in (\ref{intpse}) gives
$$
\int_0^1 \chi |u|^2 (x\del_x \pse) \>dg = O(\epsilon^\gamma).
$$
In fact, since $\psi'\ge 1$ on the interval $[\frac12,1]$ 
$$
\int_{\epsilon/2}^\epsilon \chi |u|^2 \>dg = O(\epsilon^\gamma).
$$
We conclude that $\int_0^1 \chi |u|^2\>dg$ is finite.
By Mazzeo's uniqueness result, $u=0$. 
\end{proof}

We now have all the tools in place for limiting absorption.
\begin{proposition}\label{laprop}
Fix $\zeta$ with $\re\zeta=n/2$ and choose a sequence $\zeta_j\to\zeta$ with
$\im\zeta_j= \im \zeta$ and $\re\zeta_j > n/2$.
As $\zeta_j \to \zeta$, $R_{\zeta_j}f$ converges in $x^{-\delta}L^2(X)$ to the unique
solution of 
$$
\Dz u = f
$$ 
for which $\ds u \in x^{\epsilon}L^2(X; T^*X)$.
\end{proposition}

\begin{proof}
Let $u_j = R_{\zeta_j}f$.   Suppose that $\norm{x^\delta u_j} \to \infty$.  Then we can
define a new sequence $v_j = u_j/\norm{x^\delta u_j}$ so that $\norm{x^\delta v_j} = 1$ and
$$
(\Delta - \alpha_0^2\zeta_j(n-\zeta_j))v_j = \frac{f}{\norm{x^\delta u_j}} \to 0
$$
By Proposition \ref{normdu}, $v_j$ is a bounded sequence in 
$x^{-\delta} H^1_0$.  Since the inclusion $x^{-\delta} H^1_0 \subset x^{-\delta} L^2$
is compact, any subsequence of $v_j$ has a subsequence converging in $x^{-\delta} L^2$.
The limit $v$ of and such subsequence satisfies $\Dz v = 0$ and by Proposition
\ref{radiation} it also satisfies
$x^{-\epsilon}\chi \ds v\in L^2$.  Thus by Proposition \ref{unique} $v=0$.  Since
$\norm{x^\delta v_j} = 1$ this is not possible.

The conclusion is that $u_j$ is a bounded sequence in $x^{-\delta} L^2$.  By the same
reasoning we applied to $v_j$ above, any subsequence of $u_j$ has a subsequence converging in
$x^{-\delta} L^2$ to some $u$ such that $\Dz u = f$ and $x^{-\epsilon}\chi
\ds v\in L^2$.  Since such a $u$ is unique by Proposition \ref{unique} and therefore
independent of the subsequence, we have $u_j\to u$.
\end{proof}

Using the local parametrix we can refine the estimates on $u$ considerably.
In order to avoid the singularities of $\sigma$ we introduce:
\begin{equation}\label{sigchi}\begin{split}
\chi \in\cinf(X)&:\quad \chi=0\text{ on a neighborhood of }\{\alpha^2 = \tfrac{4\lambda}{n^2}\},\\
\sigma_{reg} \in \cinf(\bX)&: \quad \sigma_{reg} = \sigma \text{ on }\supp(\chi)\cap\bX.
\end{split}\end{equation}

\begin{proposition}\label{regularity}
Suppose that $u \in x^{-\delta} L^2$ for all $\delta>0$ and $\Dz u
\in \dot\cinf(X)$.  If $u$ also satisfies the radiation condition
$(x\del_x - \sigma)u \in x^{\epsilon} L^2$ for $\epsilon>0$, then
$$
\chi u(x,y) \in \calA_{\sigma_{reg}}(X),
$$
where $\chi$, $\sigma_{reg}$ satisfy (\ref{sigchi}).
\end{proposition}

\begin{proof}
Let $P_\zeta = \Dz$ with locally defined right parametrix $M_\zeta$ as in Proposition \ref{param}.
For notational convenience we'll consider only the global case, since the introduction of 
cutoff functions for the local case is a simple matter.

Assume $u$ is as stated.
The idea is to apply the transpose of the equation $P_\zeta M_\zeta = 1 - E_\zeta$ to $u$.  To
justify this requires integration by parts, and this is where the radiation condition comes in.   Consider
the real pairing
$$
\int (P_\zeta M_\zeta \phi) u \>dg
$$
for arbitrary $\phi \in x^\epsilon L^2$.  If the integral is cutoff at $x=s$, then we can transfer
$P_\zeta$ to $u$, picking up a boundary term:
\begin{equation}\label{pmove}
\int_{x\ge s} P_\zeta (M_\zeta\phi) u \>dg = \int_{x\ge s} (M_\zeta\phi) P_\zeta u \>dg
- \int_{x=s} [M_\zeta\phi \cdot x\del_x u - u\cdot x\del_x (M_\zeta\phi)] \alpha x^{-n} dh.
\end{equation}
Let $F(s)$ denote the boundary correction from the right-hand side of (\ref{pmove}).  
The radiation condition tells us that $(x\del_x - \sigma) u \in x^\epsilon L^2$, 
and by the structure of $M_\zeta$ and Lemma \ref{maplemma} we have also $(x\del_x - \sigma) M_\zeta\phi\in
x^\epsilon L^2$.  Thus the expression in brackets in $F$ lies in $x^{\epsilon-\delta} L^1(dg)$,
because the leading terms cancel each other.  In other words
$$
\int_0 s^{\gamma-\epsilon} |F(s)| \frac{ds}s <\infty,
$$
which implies $\lim_{s\to 0} F(s) = 0$ because we can assume $\gamma<\epsilon$.  

To transfer $M_\zeta \phi$ to $M_\zeta^t u$, requires only Fubini's theorem, and so we have
$$
\int P_\zeta (M_\zeta\phi) u \>dg = \int \phi M_\zeta^t P_\zeta u \>dg.
$$
Substituting $P_\zeta M_\zeta = I - E_\zeta$ in the left-hand integral then gives
$$
\int \phi (I - E_\zeta^t) u \>dg =  \int \phi M_\zeta^t P_\zeta u \>dg.
$$
Since these properties hold for any $\phi \in x^\epsilon L^2$, we conclude that
$$
u= M_\zeta^t P_\zeta u+ E_\zeta^t u.
$$
Note that $P_\zeta u\in \dot\cinf(X)$, so $M_\zeta^t P_\zeta u \in \calA_\sigma(X)$ by 
Lemma \ref{maplemma}.  Also $E_\zeta^t \in \Psi_{\sigma_l,\infty}(X)$ so we can show $E_\zeta^t u\in
\calA_\sigma(X)$ by moving derivatives under the integral.  
Hence $u\in \calA_\sigma(X)$ as claimed. 
\end{proof}

Combine Propositions \ref{laprop} and \ref{regularity} gives the full limiting absorption
result:
\begin{theorem}\label{lathm}
Let $\re\zeta=n/2$ and $f\in\dot\cinf(X)$.  The limiting absorption principle defines a unique function
$R_\zeta f$ solving
$$
\Dz R_\zeta f = f,
$$
such that $\chi R_\zeta f \in \calA_{\sigma_{reg}}(X)$ for any $\chi$, $\sigma_{reg}$ as in (\ref{sigchi}).
\end{theorem}

\section{Scattering matrix for the irregular continuous spectrum}\label{poisson}

Fix a $\zeta$ such that $\lambda = \alpha_0^2\zeta(n-\zeta) 
\in (\frac{\alpha_0^2n^2}4, \frac{\alpha_1^2 n^2}4]$.   Given a function
$f\in\dot\cinf(W_\lambda)$, by Corollary \ref{cdsolve} we can find $u_1$
such that $(\Dl) u_1 = \phi\in \dot\cinf(X)$ and $\chi u_1\in \calA_{n-\sigma_{reg}}(X)$ with 
$x^{n-\sigma} f(y)$ as the leading term in the boundary expansion near $W_\lambda$.  (Here
$\chi,\sigma_{reg}$ are defined as in (\ref{sigchi}).)  
Then $R_\zeta \phi$ is defined by Theorem \ref{lathm} and 
$u = u_1 - R_\zeta \phi$ satisfies $(\Dl) u =0$.  We thus have proven:
\begin{proposition}\label{nuparam}
For $\zeta$ as above, given $f\in\dot\cinf(W_\lambda)$, there is a unique solution $u$ of 
$$
(\Dl) u = 0,
$$
with $\chi u \in \calA_{n-\sigma_{reg}}(X) + \calA_{\sigma_{reg}}(X)$ and 
$$
u \sim x^{n-\sigma} f + x^\sigma f' \text{ near }W_\lambda,
$$
where $f'\in\cinf(W_\lambda)$.  Thus the map $f\mapsto f'$ defines a
scattering matrix
$$
S_\zeta: \dot\cinf(W_\lambda) \to \cinf(W_\lambda).
$$
\end{proposition}
This result shows also that outside $W_\lambda$ there is no scattering.  On $\supp(\chi)\bck
W_\lambda$ we have $\re\sigma_{reg} > n/2$, and $\im\sigma_{reg} = 0$, for any $\chi,\sigma_{reg}$.
So the generalized eigenfunctions are $L^2$ near $\bX\bck \overline{W_\lambda}$, and there is no
incoming/outgoing distinction to be made.

The limiting absorption principle itself tells nothing of the structure of $R_\zeta$ or $S_\zeta$.  
In this section we will combine the local parametrix construction
with limiting absorption to obtain more information on the Poisson and scattering kernels.   As one would
expect, the local structure of these kernels near $W_\lambda$ is the same as in the case
$\zeta\in\bbC\bck\Gamma$ studied in \S\ref{merosec}.   

Without any restrictions on the set $\{\alpha^2 = \tfrac{4\lambda}{n^2}\}$, 
one can't really hope to resolve the singularities of $\sigma$ and construct the full resolvent.  
In this section we'll continue to work in full generality and just focus on the scattering kernel.
The idea will be to mimic the construction of the Poisson kernel in \S\ref{merosec},
but starting from the local parametrix instead of the resolvent.  

By Proposition \ref{param}, we have a local parametrix $M_\zeta = A_\zeta+B_\zeta$ solving the equation
$$
(\Dl) (A_\zeta + B_\zeta) = 1 - E_\zeta,
$$ 
in a neighborhood of $W_\lambda \times W_\lambda$, 
where 
\begin{equation*}\begin{split}
A_\zeta &\in \psop^{-2}(X),\\
B_\zeta &\in \psop_{\sigma_l,\sigma_r}(X),\\
E_\zeta &\in \Psi_{\infty,\sigma_r}(X).
\end{split}\end{equation*}

Let us state this precisely using cutoffs.  Fix $G\subset X$ with $G\cap \bX \subset\subset W_\lambda$ 
and let $\chi_1$, $\chi_2$ be smooth functions whose supports intersect $\bX$ inside $W_\lambda$ such that
$\chi_1 = \chi_2 = 1$ on $G$. 
If we let these act as multiplication operators, we have
\begin{equation}\label{locpar}
(\Dl) \chi_1(A_\zeta +  B_\zeta)\chi_2  = \chi_1(1 - E_\zeta)\chi_2 + [\Delta, \chi_1]
A_\zeta \chi_2 + [\Delta, \chi_1] B_\zeta \chi_2.
\end{equation}
By making the support of one cutoff fit inside the other, we are free to assume that 
$\chi_2 \cdot d\chi_1 = 0$.  Recall that by construction $A_\zeta$ can be supported in an arbitrarily narrow
neighborhood of the (lifted) diagonal.  So by manipulating these assumptions, we can ensure that 
$$
[\Delta, \chi_1] A_\zeta \chi_2=0.
$$

As in \S\ref{merosec}, we identify operators with their kernels using the Riemannian density.
The restriction map used to obtain the Poisson kernel there was
$$
r:  k(x,y,x',y') \mapsto
{x'}^{-\sigma(y')} k(x,y,x',y') |_{x'=0},
$$
which we now want to apply to (\ref{locpar}).
Actually, we'll want to consider the restriction in terms of the lifts of distributions to $\sxx$ and
$\sxbx$.  

$A_\zeta$ vanishes to infinite order at the right face, so $r(A_\zeta)=0$.  
The kernel of the multiplication operator
$\chi_1\chi_2$ lives on the lifted diagonal, and so also $r(\chi_1\chi_2)=0$.  For the residual term we
immediately conclude
\begin{equation}\label{rcz}
r(\chi_1 E_\zeta \chi_2) \in \dot\cinf(X\times \bX)
\end{equation}
Note that because $\tsl$ is an extension of the indicial root off the left face, we have
$$
(\Dl) \chi_1 B_\zeta \chi_2\text{ and } \chi_1 (\Dl)  B_\zeta \chi_2\in \calA_{\sigma_l|1,\sigma_r,0},
$$
implying also
$$
[\Delta, \chi_1] B_\zeta \chi_2 \in \calA_{\sigma_l|1,\sigma_r,0}.
$$ 
By the choice of supports, $[\Delta, \chi_1] \chi_2 = 0$, so $[\Delta, \chi_1] B_\zeta \chi_2$ is cutoff near
the diagonal in $\xx$, which makes the behavior of $B_\zeta$ at the front face irrelevant.  Hence
\begin{equation}\label{rbz}
r([\Delta, \chi_1] B_\zeta \chi_2) \in \calA_{\sigma|1}(\xbx).
\end{equation}
Here we should to use $\sigma_{reg}$ as in (\ref{sigchi}), but we'll drop the subscript for the rest of this
argument.

Finally we set
$$
M_\zeta = r(\chi_1 B_\zeta \chi_2).
$$
which is the proposed parametrix of the Poisson kernel.  Using (\ref{rcz}) and (\ref{rbz}) in (\ref{locpar})
gives the result
\begin{equation}
(\Dl) M_\zeta = F_\zeta \in \calA_{\sigma|1}(\xbx).
\end{equation}
Note that the support of $F_\zeta$ meets the boundary of $\xbx$ only in $W_\lambda\times W_\lambda$.
The error term can be improved by solving away the boundary expansion of $F_\zeta$ as in Lemma
\ref{solveaway} (with $y'$ as an extra parameter).  We thereby find $M_\zeta'\in \calA_{\sigma|1}(\xbx)$
such that
$$
(\Dl) M_\zeta' - F_\zeta = F_\zeta' \in \dot\cinf(\xbx).
$$

Finally, Theorem \ref{lathm} can be applied to remove the error term $F_\zeta'$, again with $y'$
as a parameter.  This gives us
$$
M_\zeta'' = - R_\zeta F_\zeta',
$$
with $\chi_1 M_\zeta'' \chi_2 \in \calA_\sigma(\xbx)$.
Then, setting $E_\zeta = M_\zeta -M_\zeta' + M_\zeta''$,
we have 
\begin{equation}
(\Delta + \lambda) E_\zeta = 0.
\end{equation}
Note that the structure of $E_\zeta$ is not determined near regions of the boundary where
$\alpha^2 = \frac{n^2\lambda}4$.

The action of $E_\zeta$ on a function $f\in \dot\cinf(W_\lambda)$ is given by 
$$
E_\zeta f(x,y) = \int_{W_\lambda} E_\zeta(x,y,y') f(y') \>dh_0(y'). 
$$
This clearly yields a function smooth in the interior of $X$.  To see that it is the Poisson kernel,
we need only study its asymptotic expansion at the boundary near $W_\lambda$.  
However, if we apply a cutoff then $\chi_1 E_\zeta f$ has the same structure as in the global case
considered in \S\ref{merosec}.  By the argument used in that section we have:
\begin{theorem}\label{nusmatrix}
For $\lambda = \alpha_0^2 \zeta(n-\zeta) \in (\frac{\alpha_0^2 n^2}4, \frac{\alpha_1^2 n^2}4]$,
the the Poisson kernel on $W_\lambda$ (with respect to the Riemannian density of $h_0$) is
$$
\frac{1}{B(\sigma) \sqrt{h_0}} E_\zeta.
$$
The scattering matrix $S_\zeta: \dot\cinf(W_\lambda) \to \cinf(W_\lambda)$ 
is a pseudodifferential operator with principal symbol 
$$
2^{n-2\sigma} \frac{\Gamma(\frac{n}2-\sigma)}{\Gamma(\sigma-\frac{n}2)} |\xi|^{2\sigma-n}_{h_0}
$$
\end{theorem}

\section{Resolvent for the irregular continuous spectrum}\label{resolvent}

We now turn to the question of the structure of the resolvent in the irregular part of the continuous
spectrum. In this section we will fix $\zeta$ such that $\lambda = \alpha_0^2\zeta(n-\zeta) \in
(\frac{\alpha_0^2 n^2}4, \frac{\alpha_1^2 n^2}4)$ 
and undertake a direct construction of the resolvent. 
Since analytic  continuation is not available, we will still rely on the limiting absorption principle to
prove  existence of an inverse $R_\zeta = (\Dl)^{-1}$ on the appropriate space.  The structure of the
$R_\zeta$ will be revealed by finding a sufficiently good parametrix.
Construction the parametrix requires resolution of the singularities of
$\sigma$, which is possible only if the singular set is sufficiently nice.  Henceforth
we make the assumption:
\begin{equation}\label{crassume}
\frac{2\sqrt{\lambda}}n \text{ is a regular value of }\alpha,
\end{equation}
which of course holds for generic $\lambda$ by Sard's Theorem.
Under this assumption the \textit{crossover region},
$$
\Lambda = \{y:\; \alpha^2 = \tfrac{n^2\lambda}4\} = \{y: \sigma = \tfrac{n}2\} \subset \bX
$$
(i.e. the singular set of $\sigma$), is a regular hypersurface in $\bX$.  

Away from $\Lambda$ we already have a parametrix, so we will only be concerned with local coordinates
covering neighborhoods of $\Lambda$.
For notational convenience we always choose the coordinate
$y_n$ to be a particular defining function of $\Lambda$: 
\begin{equation}\label{ynconv}
y_n = \frac{\lambda}{\alpha^2} - \frac{n^2}4.
\end{equation}
Within the coordinate neighborhood the indicial root depends solely on $y_n$,
$$
\sigma(y) = \frac{n}2 + \sqrt{-y_n}.
$$
As in \S\ref{lasec}, we introduce an extension $\tsigma$ which locally has the form
\begin{equation}\label{tsigdef}
\tsigma(x,y) = \frac{n}2 + \sqrt{ix - y_n},
\end{equation}
near $y_n = 0$.  The $ix$ term can be cutoff outside of a neighborhood of $\Lambda$ in some
arbitrary way that we won't bother to notate.

\subsection{Crossover blow-up}\label{cblup}
The singularity of $\tsigma$ is resolved by the blow-up 
$$
\Xl = [X; \Lambda].
$$  
To form $\Xl$ from $X$, $\Lambda$ is replaced by the \textit{crossover face} 
$$
\tL = SN(\Lambda) \simeq \Lambda \times S^1_+.
$$
Under our convention for local coordinates (\ref{ynconv}), the blow-up 
just amounts to the introduction of polar coordinates centered at the origin in the $(x,y_n)$ half-plane.  
Since $\bX\bck\Lambda$ is the set on which $\sigma$ is smooth, we will refer to the lift of this
set to $\Xl$ as the \textit{regular face}.
The lift of $\tsigma$ to $\Xl$ will be denoted again by $\tsigma$.

Projective coordinates show how the blow up resolves the singularity of $\tsigma$. 
In the interior of the crossover face we can use coordinates $x, w = y_n/x$.  Thus
$$
\tsigma = \frac{n}2 + x^{1/2} \sqrt{i - w},
$$
which is homogeneous in $x$ with a coefficient smooth in the interior of $\tL$.  To handle the corners on
either side of the face we use the coordinates $y_n, s = \pm x/y_n$.  For $y_n>0$, for example, we have
$$
\tsigma = \frac{n}2 + y_n^{1/2} \sqrt{is - 1}.
$$
This is smooth up to the regular face $s=0$, and again homogeneous at the crossover face.

To define the space of polyhomogeneous conormal functions
of order $\sigma$ we can't just appeal to Definition \ref{phgvo},
because $\sigma$ is not a smooth function on the regular face.  (It takes the constant value
$\frac{n}2$ on the crossover face, so no problem there.)
Fortunately, the square-root singularity in $\sigma$ is comparatively mild, and the only change
needed in the definition from \S\ref{basymsec} is to make the asymptotic expansion at the
crossover face decrease by half-integer powers rather than integer.   

Let us choose boundary defining functions $\rho_1,\rho_2$ for $\Xl$ and radial vector fields $V_1,V_2\in
\calV_b(\Xl)$, where the faces are ordered regular, crossover.  
\begin{definition}\label{phgso}
$\calA_{\sigma,[\frac{n}2]}(\Xl)$ consists of $u\in \cinf(X^\circ)$ such that
for any $m_1,m_2\in \bbN_0$ we have
$$
\Bigl[\prod_{k=0}^{m_1-1} (V_1 - k)^{k+1} \prod_{l=0}^{m_2-1}(V_2-l/2)^{\kappa_l+1}\Bigr] (x^{-\tsigma}u)
\in \calA^{(m_1,m_2/2)-}(\Xl),
$$
for some sequence $\kappa_l\in\bbN_0$ with $0=\kappa_0<\kappa_1<\dots$.
\end{definition}
The brackets in the notation $\calA_{\sigma,[\frac{n}2]}$ are meant to reflect both the half-step expansion
and the extra growth of logarithmic powers represented by $\{\kappa_l\}$. 
Applying the same reasoning as in Proposition \ref{bexp}, we see that $u$ has asymptotic expansion
as usual at the regular face, but at the crossover face powers grow by half-steps:
\begin{equation}\label{crossexp}
x^{-\tsigma} u \sim \sum_{l=0}^\infty \sum_{j\le \kappa_l} \rho_2^{l/2} (\log \rho_2)^j a_{l,j} 
\text{ near }\tL.
\end{equation}
As in \S\ref{basymsec} we define spaces with truncated expansions at the regular face by
$$
\calA_{\sigma|k,[\frac{n}2]}(\Xl) = (\rho_1\log\rho_1)^k \calA_{\sigma,[\frac{n}2]}(\Xl).
$$
For the crossover face it will prove more convenient to use the notation $[\frac{n+k}2]_+$,
where $[a]_+$ indicates that $\rho_2^{-a} u$ has an expansion of the form (\ref{crossexp}), but with
$\kappa_0>0$.

The extension $\tsigma$ depends of course on the choice of $x$, but the notation
$\calA_{\sigma,[\frac{n}2]}(\Xl)$ is meant to reflect the fact that the space is
independent of this choice. 
\begin{lemma}
The space $\calA_{\sigma,[\frac{n}2]}(\Xl)$ does not depend on the choice of coordinate $x$ 
used to define the extension $\tsigma$ in (\ref{tsigdef}).
\end{lemma}
\begin{proof}
A different choice of boundary defining function could always be written $x\phi$, where $\phi\in\cinf(X)$
and is strictly positive.  Thus we are concerned with the difference
$$
\sqrt{ix\phi - y_n} - \sqrt{ix-y_n} = \rho_1 \sqrt{\rho_2} f,
$$
for $f\in\cinf(\Xl)$.  For simplicity we can assume $\rho_1\rho_2=x$, so that we want to show that
$$
(\rho_1\rho_2)^{(\rho_1 \sqrt{\rho_2} f)}\cdot\calA_{\sigma,[\frac{n}2]}(\Xl) \subset 
\calA_{\sigma,[\frac{n}2]}(\Xl). 
$$

We use an induction similar to the proof of Lemma \ref{aindep}.  Defining the commutators
$$
Z_{q_1,q_2} = \Bigl[\prod_{k=0}^{q_1-1} (V_1 - k)^{k+1} \prod_{l=0}^{q_2-1}(V_2-l/2)^{\kappa_l+1}, \;
(\rho_1\rho_2)^{(\rho_1 \sqrt{\rho_2} f)} \Bigr],
$$
the result will follow if we can show
$$
Z_{q_1,q_2}: \calA_{\sigma,[\frac{n}2]}(\Xl) \to \calA_{\sigma|q_1,[(n+q_2)/2]_+}(\Xl).
$$
The operators $Z_{1,0}$ and $Z_{0,1}$ are functions that act by multiplication,
and a simple computation shows that
$$
Z_{1,0} \in (\rho_1\log\rho_1) \sqrt{\rho_2} (\rho_1\rho_2)^{(\rho_1 \sqrt{\rho_2} f)} \cdot
\cinf(\Xl),
$$
and
$$
Z_{0,1} \in (\sqrt{\rho_2}\log\rho_2) \rho_1 (\rho_1\rho_2)^{(\rho_1 \sqrt{\rho_2} f)} \cdot
\cinf(\Xl),
$$
so the mapping result holds for these cases by the characterization in terms of asymptotic expansions.
The inductive step follows precisely as in Lemma \ref{aindep}, once we note that
\begin{equation*}\begin{split}
Z_{1,0}: \calA_{\sigma|q_1,[(n+q_2)/2]_+}(\Xl) &\to \calA_{\sigma|(q_1+1),[(n+q_2)/2]_+}(\Xl) \\
Z_{0,1}: \calA_{\sigma|q_1,[(n+q_2)/2]_+}(\Xl) &\to \calA_{\sigma|q_1,[(n+q_2+1)/2]_+}(\Xl) 
\end{split}\end{equation*}
\end{proof}

\begin{lemma}\label{vcala}
The space $\calA_{\sigma,[\frac{n}2]}(\Xl)$ is preserved under the action of lifts of 
vector fields from $\calV_0(X)$ to $\Xl$.
\end{lemma}
\begin{proof}
It is easy to verify that vector fields in $\calV_0(X)$ lift to vector fields in $\calV_b(\Xl)$
by considering local coordinates.  For example, in the interior of the crossover face the coordinates 
$x, w = y_n/x$ are valid.   The lifts of the vector fields $x\del_x$ and $x\del_{y_n}$ to $\Xl$ are
\begin{equation*}\begin{split}
x\del_x &\longrightarrow x\del_x - w\del_w\\
x\del_{y_n} &\longrightarrow \del_w
\end{split}\end{equation*}
At the boundary of the regular and crossover faces, we use coordinates $y_n$, $t = x/y_n$,
and the lifts are
\begin{equation*}\begin{split}
x\del_x &\longrightarrow t\del_t\\
x\del_{y_n} &\longrightarrow y_nt\del_{y_n} - t^2\del_t
\end{split}\end{equation*}
The lemma follows immediately from the invariance of $\calA_{\sigma,[\frac{n}2]}(\Xl)$ under
$\calV_b(\Xl)$.
\end{proof}

The proof of Lemma \ref{vcala} reveals one of the complications of introducing $\Xl$.  The
lifts of vector fields from $\calV_0(X)$ vanish at the regular face, but are only tangent
to the crossover face.  Thus the indicial equation, which defined $\sigma$ in the first place,
is not valid on $\Xl$.  To put this another way, the highest order terms in $(\Dl)x^\tsigma f$
cancel out if $f \in \cinf(X)$ but not for $f\in\cinf(\Xl)$.
Thus the analog of Corollary \ref{cdsolve} in this context is:
\begin{proposition}\label{xlsolve}
Given $f \in \cinf(\bX)$ can find $u\in\calA_{\sigma,[\frac{n}2]}(\Xl)$ satisfying
$$
(\Dl) u \in \calA_{\infty, [\frac{n+1}2]_+}(\Xl)
$$
(for $\zeta$ such that $\lambda$ satisfies (\ref{crassume})),
so that $x^{\tsigma} f$ gives the leading term in the asymptotic expansions of $u$ at both the regular or
crossover faces.
\end{proposition}
\begin{proof}
Extending $f$ smoothly into the interior in some arbitrary way, we have
$$
(\Dl) x^{\tsigma} f  \in \calA_{\sigma|1, [\frac{n+1}2]_+}(\Xl),
$$
because $\sigma$ is the indicial root.  As in the proof of Lemma \ref{vcala}, let us use coordinates
$y_n$, $t = x/y_n$.  By the lifts of vector fields computed there we see that
\begin{equation*}\begin{split}
(\Dl) t^{\tsigma+k} (\ln t)^l a_{k,l} &= \alpha^2 k (n-2\tsigma-k) t^{\tsigma+k} (\ln t)^l a_{k,l}\\
&\qquad + t^{\tsigma+k} (\ln t)^l \Bigl[\sqrt{y_n} (t\ln t) g_1 + tg_2 + y_n (t\ln t)^2 g_3  \Bigr],
\end{split}\end{equation*}
where $g_1, g_2, g_3 \in \cinf(\Xl)$.  This shows that the asymptotic expansion of $(\Dl) x^{\tsigma} f$ 
at the regular face may be solved away exactly as in Lemma \ref{solveaway}, giving $u$ as stated.
\end{proof}

\subsection{Normal operator}
An essential tool in the Mazzeo-Melrose construction of the resolvent \cite{MM}
is the normal operator defined by restricting a kernel to the front face of $\sxx$. 
The definition is most conveniently made using half-densities.  In particular they
allow us to easily write local formulas for invariant expressions.  
The Riemannian half-density on $X$ determined by
$g$ is a smooth section of the singular density bundle
$\shlf = x^{-(n+1)/2}\Omega^{1/2}$, where $\Omega$ is the density bundle.  By the same convention we 
define $\shlf(\xx) = (xx')^{-(n+1)/2}\Omega$.  We'll also let $\shlf(\sxx)$ denote
the lift of $\shlf(\xx)$ to the stretched product.  The Laplacian on half-densities is defined by
$$
\hD(f\mu_g^{1/2}) = (\Delta f) \mu_g^{1/2},
$$
where $\mu_g$ is the Riemannian density.  We'll work with half-densities for the rest of this section. 
Operators will be assumed to act on $\cinf(X; \shlf)$, 
so their kernels are naturally interpreted as distributional half-densities in $C^{-\infty}(\xx; \shlf)$.

Let $F$ be the front face of $\sxx$, and for $q\in \bX$
let $F_q$ be the fiber of $F$ over $q$. If $K$ is an operator in 
$$
\psop^{m}_{\sigma_l,\sigma_r}(X; \shlf) = \psop^{m}(X; \shlf) + \psop_{\sigma_l,\sigma_r}(X; \shlf),
$$ 
then the normal operator of $K$ at $q$ is
$$
N_q(K) = \beta^* K|_{F_q} \in \psop^m_{\sigma(q), \sigma(q)}(F_q; \shlf(\sxx)|_{F_q}).
$$
The key fact here is that the $F_q$ can be identified with the group of boundary preserving linear
transformations of $X_q = T^+_q X$, while the restriction of $\shlf(\sxx)$ to $F_q$ is naturally identified
with $\shlf(X_p\times X_p)$.  Thus $N_q(K)$ acts naturally as a convolution operator 
$$
N_q(K): \dot\cinf(X_q; \shlf) \to  C^{-\infty}(X_q; \shlf).
$$
To write this map down concretely, let $(x,y,x',y')$ be the coordinates on $\xx$.
For $\sxx$ we use $x,y$ and the projective coordinates $t = x/x'$, $u=(y-y')/x'$.  Fixing a base point 
$q = (0,y_0) \in \bX$,
$(t,u)$ become coordinates for the fiber $F_q$.  It is convenient to give $X_q$ the coordinates $x,y$ given
by linearizing the functions $x$ and $y$.  The action of $F_q$ on $X_q$ can then be written
$$
(x,y) \cdot (t,u)^{-1} = (\frac{x}{t}, y - \frac{x}{t} u).
$$
We'll introduce reference half-forms $\mu = |\frac{dx\>dy}{x^{n+1}}|^{1/2}$, as a section of
$\shlf(X_q)$, and 
$$
\gamma = \Bigl| \frac{dx\>dy\>dt\>du}{x^{n+1} t}\Bigr|^{1/2} \in \cinf(\sxx; \shlf).
$$
Let $\beta^* K$ be given locally by $k(x,y,t,u)\cdot\gamma$.  Then the convolution action of the normal operator
is
$$
N_q(K) (f\cdot\mu) = \int k(0,y_0,t,u) f(\frac{x}{t}, y - \frac{x}{t} u) \>\frac{dt\>du}{t} \cdot \mu,
$$
where $f\in\cinf(X_q)$.  

The usefulness of this definition rests on the following facts. 
Let $\Diff^m_0(X)$ be the space of differential operators on $X$ generated by $\calV_0$.  
If $\Diff^m_0(X; \shlf)$ is the corresponding space of operators on half-densities, 
then we can naturally regard
$\Diff^m_0(X; \shlf)$ as a subspace of $\psop^m_{\infty,\infty}(\sxx; \shlf)$.

\begin{lemma}\label{npk} \cite{MM}
For $P\in \Diff^m_0(X; \shlf)$ the normal operator $N_q(P)$ is given by ``freezing coefficients'' at $q$. 
That is, for 
$$
P(f\cdot \mu) = \sum_{k+|\alpha|\le m} p_{k,\alpha}(x,y) (x\del x)^k (x\del_y)^\alpha f \cdot \mu,
$$
we have (abusing notation by using the same coordinates for $X$ and $X_q$):
$$
N_q(P)(f\cdot \mu) = \sum_{k+|\alpha|\le m} p_{k,\alpha}(0,y_0) (x\del x)^k (x\del_y)^\alpha f \cdot \mu. 
$$
Moreover, if $K\in\psop^\infty_{a,b}(\sxx; \shlf)$, 
\begin{equation}\label{npkline}
N_q(PK) = N_q(P) \cdot N_q(K).
\end{equation}
\end{lemma}

\subsection{Crossover stretched product}
To construct the resolvent when $\sigma$ has singularities, we must combine the
stretched product with the blow-up of $\Lambda$ in \S\ref{cblup}.  To understand what is needed here, we
recall some facts from the general theory developed in Chapter 5 of \cite{MeDA}.
A \textit{p-submanifold} $Z$ of a manifold with corners $M$ is a submanifold such that
near each point of $Z$ there are coordinates $(x,y) \in \bbR_+^k\times \bbR^{m-k}$ 
for $M$ such that locally
$$
Z = \{x_{l+1}=\dots=x_k=0,\; y_{k+j+l} = \dots = y_m =0\}
$$
(`p' stands for product).  All of the submanifolds we want to blow up will be p-submanifolds.
If $Z$ and $Y$ are p-submanifolds of $M$ then the lift of $Z$ to $[M;Y]$ through the blow-down map
$\beta: [M; Y]\to M$ is defined in two distinct cases:
$$
\beta^* Z = 
\begin{cases}
\beta^{-1}Z& \text{if }Z\subset Y,\\
\overline{\beta^{-1}(Z\bck Y)}& \text{if }Z = \overline{Z\bck Y}.
\end{cases}
$$
Assuming that the lift $\beta^*Z$ is a p-submanifold of $[M; Y]$, we can define
the double blow-up $[M; Y; Z]$ as $[[M;Y]; \beta^*Z]$.

\begin{lemma}\label{buorder}\cite{MeDA}
The double blow-ups $[M;Y;Z]$ and $[M;Z;Y]$ are well-defined and naturally diffeomorphic by the unique
extension of the identity map on the interior under the following conditions:
\begin{enumerate}
\item $Y$ is a p-submanifold of $M$ and $Z$ is a p-submanifold of $Y$.
\item $Y$ and $Z$ are p-submanifolds of $M$ meeting transversally.
\end{enumerate}
\end{lemma}

\begin{figure}
\centerline{\epsfysize=2.3in \epsfbox{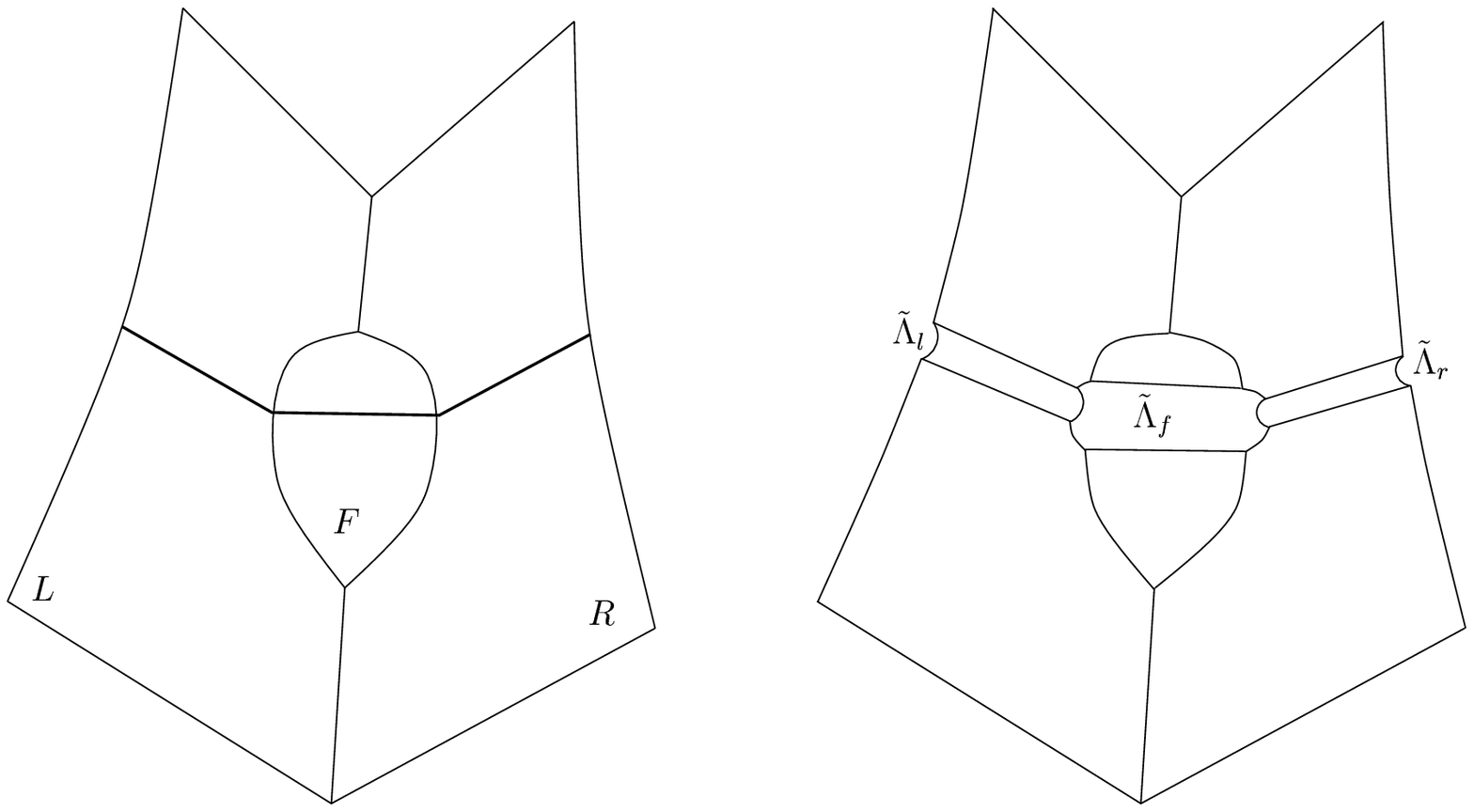}}
\caption{Copies of the crossover set in $\sxx$ are blown up to form $\esxx$.}
\label{esxx}
\end{figure}

In $\xx$ the submanifolds needing to be resolved are $S = \{x=x'=y-y'=0\}$, 
$\Lambda_l = \{x=y_n=0\}$, and $\Lambda_r = \{x'=y_n'=0\}$.  In order to make
a symmetric definition we first need to blow-up the intersection 
$$
\Lambda_f = \Lambda_l \cap S = \Lambda_r\cap S = \{x=x'=y-y'=y_n=0\}.
$$
The \textit{crossover stretched product} is 
\begin{equation}\label{esxxdef}
\esxx = [\xx; \Lambda_f; S; \Lambda_l; \Lambda_r].
\end{equation}
There are many equivalent ways to order the blow-ups, though not all possibilities are permitted. 
For example, since $\Lambda_f \subset S$, we can interchange the first two blow-ups and write
$\esxx = [\sxx; \Lambda_f; \Lambda_l; \Lambda_r]$.   This is illustrated in Figure \ref{esxx}.
The crossover stretched product has six faces.  The left and right faces of $\xx$ lift to \textit{left} and 
\textit{right regular faces}, and the lift of $S$ is the \textit{front regular face}.
The lifts of $\Lambda_f, \Lambda_l, \Lambda_r$ are the \textit{front, left} and \textit{right
crossover faces}, respectively.

The lift of $S$ to $[\xx; \Lambda_f]$ does not intersect the lifts of $\Lambda_l$
and $\Lambda_r$, and the intersection of $\Lambda_l$ with $\Lambda_r$ is clearly transverse.  
Thus in (\ref{esxxdef}) the last three blow-ups may be taken in any order.   
Because of this there are natural projections $\esxx \to \Xl$ through the right and left
factors.  To project on the left, for example, we can first blow down $\Lambda_r$ and 
$S$ to get $[\xx; \Lambda_f; \Lambda_l]$.  
Then because $\Lambda_f\subset \Lambda_l$ we can interchange these two blow-ups
and blow down $\Lambda_f$, leaving $[\xx; \Lambda_l] = \Xl \times X$.  
Finally we project onto $\Xl$ to give the full map.  

To describe expansions at the crossover faces we will continue to use the notations $[a]$ and 
$[a]_+$ introduced in \S\ref{cblup}.  The full set of boundary information for the faces will 
be given as an \textit{index family}.  For example, let
$\calM = (\sigma_l,\sigma_r,0,[\frac{n}2],[\frac{n}2],[0])$, where the order is left, right, front
regular, then left, right, front crossover.  The space $\calA_\calM(\esxx)$ is defined analogously
to $\calA_{\sigma,[\frac{n}2]}(\Xl)$ in Definition \ref{phgvo}, 
using extensions $\tsl, \tsr$ which are lifts of
$\tsigma$.  By the same arguments from \S\ref{cblup} we have:
\begin{lemma}
$\calA_\calM(\esxx)$ is well-defined independently of the choice of extensions of $\sigma_l$ and
$\sigma_r$ and is invariant under $\calV_b(\esxx)$.
\end{lemma}

\noindent
Truncated expansions at the regular faces will be notated as before, by
replacing $\sigma_l$ with $\sigma_l|k$.

We continue to use $\shlf$ to denote the half-density bundle spanned by the Riemannian
half-density on $X$ (and lifts and combinations thereof).  
For example, $\shlf(\Xl)$ denotes the lift of $\shlf(X)$ to $\Xl$.  If
$\rho, s$ are defining functions for the regular and crossover faces, respectively, then
$$
\shlf(\Xl) = \rho^{-\frac{n+1}2} s^{-\frac{n}2} \Omega^{1/2}(\Xl).
$$
The formula for $\shlf(\esxx)$ is similar.

Abusing notation slightly, we write $A\in \calA_\calM(\esxx; \shlf)$ for an operator
on sections of $\shlf(X)$, meaning that the distribution kernel of $A$ lifts to an 
element of this space.  

The front regular face still fibers over the regular face of $\Xl$, and there the definition and properties
of the normal operator are unchanged.   However, we will not be able to use this construction at
the front crossover face, beyond the first iteration.

\begin{proposition}\label{rparam}
For $\zeta$ such that $\lambda = \alpha_0^2\zeta(n-\zeta) \in (\frac{\alpha_0^2n^2}4, \frac{\alpha_1^2n^2}4)$,
there exists a parametrix $W$ for $\Dl$ such that
$$
(\hDl) W = I - F 
$$ 
where
\begin{equation*}\begin{split}
W &\in \psop^{-2}(X; \shlf) + \calA_\calM(\esxx; \shlf),\\
F &\in \calA_\calF(\esxx; \shlf),
\end{split}\end{equation*}
with $\calM = (\sigma_l,\sigma_r,0,[\frac{n}2],[\frac{n}2],[0])$ and $\calF =
(\infty,\sigma_r,\infty,[\frac{n+1}2]_+,[\frac{n}2],[\frac12]_+)$
\end{proposition}
\begin{proof}
The first part of the construction is the symbolic inversion along the lifted diagonal in $\sxx$.
This does not involve $\sigma$ and so can be done just as in \cite{MM}.  We obtain 
$A\in \psop^{-2}(X; \shlf)$ such that
$$
(\hDl) A - I = Q \in \psop^{-\infty}(X; \shlf).
$$
The kernel of $Q$, lifted to $\sxx$, lies in $\rho^\infty {\rho'}^\infty\cinf(\sxx; \shlf)$.  The next step
in the construction is to try to solve away this error term at the front face.

To this end, we seek an operator $B_0$ such that
\begin{equation}\label{ndbq}
N(\hDl) N(B_0) = N(Q).
\end{equation}
Here $N_p(Q) \in \dot\cinf(F_p)$.
For $p\in\bX$,
$$
N_p(\hDl) = \alpha^2(p) \Delta_{\beta_p} - \lambda,
$$ 
where $\beta_p$ is the metric on $X_p$ given by 
$$
\beta_p = \frac{dx^2 + h_p}{x^2},
$$
with $h_p$ the constant tensor on $X_p$ defined by the restriction of $h$ to $T^+_pX$. 
Since $h_p$ is constant, $\beta_p$ is the pull-back of the standard hyperbolic metric on $\bbH^{n+1}$
by a linear change of coordinates.  Let 
$$
R_0(\beta_p; \eta) = [\hD_{\beta_p} - \eta(n-\eta)]^{-1}
$$
denote the resolvent for this metric, which is just given by applying the linear coordinate change to 
the standard resolvent on $\bbH^{n+1}$.  
Noting that $\lambda/\alpha^2 = \sigma(n-\sigma)$, we set
$$
N_p(B_0) = \frac1{\alpha^2(p)} R_0(\beta_p; \sigma(p)) N_p(Q).
$$

For each $p$ we have $N_p(B_0) \in \calA_{\sigma(p),\sigma(p)}(F_p; \shlf|_{F_p})$.  
Although $N_p(Q)$ depended smoothly on $p$, $N_p(B_0)$ does not because of $\sigma$.   
So this is the first point at which we need the lift from $\sxx$ to $\esxx$.  
Let $\tilde S$ denote the front regular face (the lift of $S$).
Since the dependence of $N(B_0)$ on $\sigma$ is analytic, the lift of $N(B_0)$ to $\tilde S$
has a square root singularity at the crossover boundary coming from that of $\sigma$.  
Hence we have
$$
N(B_0)|_{\tilde S} \in \calA_{\sigma,\sigma,0}(\tilde S; \shlf|_{\tilde S}),
$$
where the boundaries of $\tilde S$ are left/right regular and front crossover, and
$\sigma$ is lifted to $\tilde S$ through the projection $\tilde S \to \bX\bck\Lambda$. 

The front crossover face, denoted $\tL_f$, fibers over lift of $\Lambda$
to the front face $S$ in $\sxx$.
The lift of $N(B_0)$ from $S$ to the front crossover face in $\esxx$ will be constant on such fibers.
Thus the lift of $N(B_0)$ has no singularity at the front regular boundary of $\tL$, and 
$$
N(B_0)|_{\tL_f} \in \calA_{\sigma,\sigma,0,[\frac{n}2],[\frac{n}2]}(\tL_f; \shlf|_{\tL_f}).
$$
The boundary faces of $\tL_f$ are left/right/front regular and left/right crossover
(see Figure \ref{esxx}).

By Proposition \ref{bexp}, this shows that the lift of $N(B_0)$ to $\esxx$ has the properties of 
the restriction of a kernel in $\calA_\calM(\esxx; \shlf)$ to the two front faces.  We may therefore choose
an  extension $B_0 \in \calA_\calM(\esxx; \shlf)$, having the lift of $N(B_0)$ as leading
coefficient at these faces.  Furthermore, since $N(B_0)$ was constant on fibers
of the front crossover face, we can insist that the leading coefficient of $B_0$ at
the left crossover face be constant on the fibers of this face as well.

The error term  at the next stage is
$$
Q_1 = (\hDl) B_0 - Q.
$$
By construction, the kernel of $Q_1$ vanishes at the front face, so 
$$
Q_1\in \calA_{\calM_1}(\esxx; \shlf),\qquad \calM_1 = (\sigma_l, \sigma_r, 0|1, [\tfrac{n}2],
[\tfrac{n}2],[\tfrac12]_+).
$$ 
However, there is additional decay on the front face, because $\sigma$ is
the indicial root.  The composition $(\hDl) B_0$ has kernel given by lifting
$\hDl$ to $\esxx$ through the left and applying it to the kernel $B_0$.  

Consider local coordinates $(t,z,x',\bar y, w)$, where $t = \frac{x}{y_n}$, $z = \frac{y-y'}{x'}$,
$w = \frac{y_n}{x'}$, and $\bar y = (y_1,\dots, y_{n-1})$.  This system is valid in a 
neighborhood of the intersection $\tilde L \cap \tL_f\cap \tL_l$.  The lift of vector fields
through the left is easily computed:
\begin{equation*}\begin{split}
x\del_x &\longrightarrow t\del_t \\
x\del_{\bar y} &\longrightarrow wt \del_{\bar z} + wtx' \del_{\bar y}\\
x\del_{y_n} &\longrightarrow wt\del_w - t^2 \del_t + wt \del_{\bar z}
\end{split}\end{equation*}
Here $t$ defines the left regular face, and $w$ the left crossover.  It is easy to compute that
\begin{equation}\label{lindicial}
\begin{split}
(\beta_l^*\Delta-\lambda) (tw)^{\tsigma_l+k}f &= \alpha^2 k(n-2\sigma_l-k) (tw)^{\tsigma_l+k}f\\
&\quad+(tw)^{\tsigma_l+k} \Bigl[(wx')^{1/2} (t\log tw) g_1 +  wx'(t\log tw)^2 g_2\Bigl] f\\
&\quad+(tw)^{\tsigma_l+k} tV f ,
\end{split}
\end{equation}
where $g_1,g_2 \in \cinf(\esxx)$ and $V \in \Diff_b^*(\esxx)$.  In particular, with $k=0$ we see
that $(\hDl) B_0$ vanishes at the left regular face.

At the left crossover face we have the same issue as in Proposition \ref{xlsolve}, because the final term on
the right in (\ref{lindicial}) is not of lower order in $w$.  More specifically, $tVf$ contains terms of
the form $(t\del_t)^i f$ times smooth coefficients.  However, in these coordinates $t$ is
the fiber variable for the left crossover face.  Thus choosing $B_0$ to be constant on
the fibers means that in the calculation of $(\hDl) B_0$ the term represented by $tV$ does not contribute.
We conclude that 
$$
\calM_1 = (\sigma_l|1, \sigma_r, 0|1, [\tfrac{n+1}2]_+, [\tfrac{n}2],[\tfrac12]_+).
$$

At the front regular face we may iterate the procedure of solving away terms with the normal operator,
exactly as in \cite{MM}.  
We would next apply the model resolvent to $N_p((x'\log x')^{-1}Q_1)$ (using Proposition
\ref{fibermodel} to handle the extra logarithms) and so obtain $N_p((x'\log x')^{-1} B_{1,1})$, from which
we determine
$B_{1,1}$, and so on.  (We cannot do the same at the crossover face, because the derivatives along the fiber
appearing in (\ref{lindicial}) cannot be assumed to vanish for any terms beyond $B_0$.  As a result, we
could not apply Proposition \ref{fibermodel} to the higher $Q_j$'s because of the lack of decay at the left
crossover face.)

After repeated applications of the normal operator at the front regular face we end up with a sequence
$B_{j,l}$ which can be summed asymptotically to give $B \in \calA_\calM(\esxx; \shlf)$ such that
$$
(\hDl) B - Q = G \in \calA_\calG(\esxx; \shlf),
$$
where $\calG = (\sigma_l|1, \sigma_r,\infty, [\frac{n+1}2]_+, [\frac{n}2],[\frac12]_+)$.
Finally, to get $F$ we can appeal to (\ref{lindicial}) repeatedly to solve away the power series 
at the left regular face, as in Proposition \ref{xlsolve}.
This yields $C \in \calA_\calM(\esxx; \shlf)$ such that
$$
(\hDl) C - G = F \in \calA_\calF(\esxx; \shlf),
$$
where $\calF$ is as defined in the statement.  Setting $W = A-B+C$ 
gives the result.
\end{proof}

\subsection{Compositions}\label{compss}
To refine the parametrix $W$ further, we would like to
asymptotically sum the Neumann series for $(I-F)^{-1}$ and
right multiply this by $W$.  Analyzing the result requires some composition properties 
for operators of these types.  Composition can always be broken down into a combination of pull-backs
and push-forwards, and we need to be sure that  the character of asymptotic expansions at the
boundary under such operations.

Fortunately, general results on pull-back and push-forward of polyhomogeneous conormal functions were
established in \cite{Me92, MeDA}.   We'll review briefly the results we need, which are all in the
context of constant indices. 
Although one could extend this theory to polyhomogeneous functions with variable order, complications 
would arise because of the possible crossing of indices. 
We won't undertake this here, as the constant index theory is sufficient for our purposes.

To specify a general polyhomogeneous conormal function, we give an \textit{index set} $E \subset \bbC
\times \bbN_0$ at each boundary face, such that at this face the boundary expansion has the
form
$$
\sum_{(s,k)\in E} \sum_{l=0}^k t^{a} (\log t)^l a_{s,l},
$$
where $t$ is the defining function.  To ensure a finite number of logarithmic terms at
each stage it is required that for any sequence $(s_j,k_j) \in E$ with $|(s_j,k_j)|\to\infty$
we have $\re s_j \to \infty$.  It is convenient to write
$\re E > m$ to mean  $\re a> m$ for all $(a,k)\in E$.
At this level of generality, an index family is a collection $\calE = (E_1,E_2,\dots)$, one index
set for each boundary face.  The space $\calA_\calE(M)$ is defined just
as in \S\ref{basymsec}.

On a manifold with corners $M$ one can associate to $\calV_b(M)$ a corresponding tangent bundle
$^bTM$.  A map $f:M_1\to M_2$ induces
$^bf_*: {^bT}M_1 \to {^bT}M_2$ by extension from the interior.
Let $r_j$, $\rho_j$ be defining functions for the boundary hypersurfaces of $M_1$ and $M_2$ respectively.
The map $f$ is called a \textit{b-map} if there exist $e(i,j) \in \bbN_0$ such that
$$
f^*\rho_i = h \prod r^{e(i,j)},
$$
Actually, in all of our cases $e(i,j)\in \{0,1\}$. The b-map condition is sufficient for pull-back
of polyhomogeneity.  
For $\calF$ an index family for $M_2$,  let $f^\sharp(\calF) = (E_1,E_2,\dots)$,
where 
$$
E_j = \{(b,p): b = \sum_i e(i,j) a_i,\; p = \sum_{e(i,j)\ne0}k_i,\; (a_i,k_i) \in F_i\} .
$$
\begin{proposition}
For a b-map $f: M_1\to M_2$, pull-back gives a map
$$
f^*: \calA_{\calF}(M_2) \to \calA_{f^\sharp(\calF)}(M_1) 
$$
\end{proposition}
The proof is quite straightforward: one can push forward the radial vector fields from $M_1$ to
$M_2$.

Extra conditions are needed for push-forward of polyhomogeneity.  A b-map $f:M_1\to M_2$ is called a
\textit{b-fibration} if: (1) $^bf_*$ is surjective at each point of $\del M_1$; and (2) no boundary
hypersurface of $M_1$ is mapped to a corner of $M_2$.  The latter condition means  for each $j$ there is at
most one $i$ for which $e(i,j) \ne 0$.  The b-density bundle is defined by
$\Omega_b = (\prod \rho_j)^{-1} \Omega$.  

In push-forwards the powers of logarithms may increase when two or more faces map to the same
face.  So for two index sets $E,F$ we define the \textit{extended union} 
$$
E \xU F = E\cup F \cup \{(a,k+k'+1):\; (a,k)\in E,\; (a,k') \in F\}.
$$
Then define $f_\sharp(\calE) = (F_1,F_2,\dots)$, in the special case that all 
$e(i,j) \in \{0,1\}$, by
$$
F_i = \xU_{e(i,j)\ne0} E_j.
$$
(If $e(i,j)>1$ then the orders in $E_j$ would be divided by $e(i,j)$.)

\begin{proposition}\label{pushfor}
Let $f:M_1\to M_2$ be a b-fibration and suppose that $\re E_j>0$ for all $j$ such that
$e(i,j)=0$ for each $i$ (i.e. for all $j$ such that the $j$'th boundary face maps to the interior).   
Then push-forward gives a map
$$
f_*: \calA_{\calE}(M_1; \Omega_b) \to \calA_{f_\sharp(\calE)}(M_2; \Omega_b) 
$$
\end{proposition}

In order to compose operators whose kernels live naturally on $\esxx$ we need to define a blown-up
triple product $\esxxx$ such that projections to $\esxx$ through various factors
exist and are b-fibrations.  This is analogous to the edge triple product in \cite{Ma91b}.
The philosophy behind the construction is the same as in the
definition of $\esxx$.  $X^3$ contains three copies each of $S$ and $\Lambda$, but
before blowing them up we must blow up any non-transverse intersections.
It is perhaps easiest to describe this in local coordinates, which also make evident the fact
that all submanifolds introduced are p-submanifolds.  Corresponding to the usual $(x,y)$ on $X$, we
use $(x,y,x',y',x'',y'')$ as local coordinates for $X^3$.  Let us label the three boundaries of 
$X^3$ as $L,M,R$ for left, middle, and right.  The three copies of $S$ are
\begin{equation*}\begin{split}
S_{LM} &= \{x=x'=y-y'=0\},\\ 
S_{MR} &= \{x'=x''=y'-y''=0\},\\
S_{LR} &= \{x=x''=y-y''=0\},
\end{split}\end{equation*}
and the copies of the crossover set are
$$
\Lambda_L = \{x=y_n=0\}, \qquad \Lambda_M = \{x'=y_n'=0\}, \qquad \Lambda_R = \{x''=y_n''=0\}.
$$
The intersections among the copies of $\Lambda$ are transverse, but none of the intersections involving
copies of $S$ are.  
Additional blow-ups are needed for the triple intersections,
$$
S_{LMR} = S_{LM} \cap S_{MR} \cap S_{LR}, \qquad \Lambda_{LMR} = \Lambda_L \cap S_{LMR},
$$
and double intersections
$$
\Lambda_{LM} = \Lambda_{L} \cap S_{LM}, \qquad
\Lambda_{MR} = \Lambda_{M} \cap S_{MR}, \qquad
\Lambda_{LR} = \Lambda_{L} \cap S_{LR}.
$$
The blown-up triple product is thus
$$
\esxxx = [X^3; \Lambda_{LMR}; S_{LMR}; \Lambda_{LM}; \Lambda_{MR}; \Lambda_{LR}; S_{LM}; S_{MR};
S_{LR};  \Lambda_{L}; \Lambda_{M}; \Lambda_{R}]
$$
This gives a manifold with 14 faces.  The faces will be denoted by putting a tilde over the corresponding
submanifold of $X^3$, e.g. $\tilde L$ is the lift of $\bX\times X\times X$

The claim is that these blow-ups allow the definition of three b-fibrations, 
$$
\psi_{LM}, \psi_{MR}, \psi_{LR}: \esxxx\to \esxx,
$$  
where the subscripts indicate which factors of $X$ are involved.
For example, to get $\psi_{LM}$ we blow down, in order, $\Lambda_{M}$, $\Lambda_{MR}$,
$\Lambda_{LR}$,  $S_{MR}$, $S_{LR}$, $S_{LMR}$, $\Lambda_{LMR}$, using Lemma \ref{buorder} to justify
the interchange of blow-ups from the original definition of $\esxxx$.  
This leaves us with $[X^3; \Lambda_{LM}; S_{LM}; \Lambda_L;
\Lambda_R]$ which equals $(\esxx)\times X$, and so the definition of $\psi_{LM}$ 
is completed with a projection.  To check
that these maps are b-fibrations is a simple exercise.

In order to apply the pull-back and push-forward formulas,
we need to examine the lifts of boundary defining functions from $\esxx$ to $\esxxx$.
Let $\rho_l,\rho_r,\rho_f$ be defining functions of the regular front,
right, and crossover faces of $\esxx$, and $s_l,s_r,s_f$ defining functions for the
crossover faces.  On $\esxxx$ we will denote defining functions of the regular faces by $r$ and of 
crossover faces by $\gamma$, with subscripts to indicate particular faces.  For example, $r_L$ will 
be the defining function for $L$, $r_{LM}$ for $S_{LM}$,
$r_{LMR}$ for $S_{LMR}$, etc.  And similarly $\gamma_L$ for $\Lambda_L$, $\gamma_{LM}$
for $\Lambda_{LM}$, etc.
By tracing the definitions of the $\psi$'s backwards, we can compute the lifts of the defining functions
from $\esxx$ to $\esxxx$.  The regular and crossover defining functions behave separately and completely
analogously to each other, so we can record the lifts in
a single chart:
\begin{equation}\label{liftchart}
\begin{tabular}{lccccccc}
\qquad & L\qquad& M\qquad& R\qquad& LM\qquad& MR\qquad&LR\qquad&LMR\\
$\psi_{LM}^*$&l&r&--&f&r&l&f\\
$\psi_{MR}^*$&--&l&r&l&f&r&f\\
$\psi_{LR}^*$&l&--&r&l&r&f&f\\
\end{tabular}
\end{equation}
This means, for example, that $\psi_{LM}^*(\rho_l) = r_L r_{LR}$
and also $\psi_{LM}^*(s_l) = \gamma_L \gamma_{LR}$, up to smooth non-vanishing functions which we will
be dropping throughout this discussion.  

Let $\nu_0$ be a smooth non-vanishing section of $\shlf(\esxx)$.  Then for operators represented as
kernels $A\nu_0$ and $B\nu_0$ the formula for composition is:
\begin{equation}\label{compdef}
(A\circ B)\nu_0 = \psi_{LR*} \Bigl[ \psi_{LM}^*(A\nu_0)\; \psi_{MR}^*(B\nu_0) \Bigr].
\end{equation}
In order to apply the push-forward theorem we need to relate everything to the b-density bundle $\Omega_b$. 
Let $\omega_b$ and $\omega_{b3}$ be sections of $\Omega_b(\esxx)$ and $\Omega_b(\esxxx)$ respectively.
\begin{lemma}
Written as a push-forward of b-densities, the composition formula (\ref{compdef}) becomes
\begin{equation}\label{abomega}
(A\circ B)  \omega_b = \psi_{LR*} \Bigl[ \psi_{LM}^*(A) \psi_{MR}^*(B) r_M^{-n} \gamma_M^{-n+1} \omega_{b3}
\Bigr]
\end{equation}
\end{lemma}
\begin{proof}
For convenience let us use subscripts to denote lifts by the $\psi$'s, e.g. $A_{LM} = \psi_{LM}^*A$.
Multiplying (\ref{compdef}) on each side by $\nu_0$ gives
\begin{equation}\label{abnu}
(A\circ B) \nu_0^2 = \psi_{LR*} \Bigl[ A_{LM} B_{MR} (\nu_0)_{LM}
(\nu_0)_{MR} (\nu_0)_{LR} \Bigr]
\end{equation}
If $\mu$ is a smooth half-density on $\xx$, then 
\begin{equation}\label{mu3}
\mu_{LM}\mu_{MR}\mu_{LR} = \mu_3^2,
\end{equation}
where $\mu_3$ is a smooth half-density on $\esxxx$.   Let $\beta:\esxx\to\xx$ be the blow-down.  By
computing the codimensions of the submanifolds to be blown-up, we can see that
$$
\beta^*\Omega^{1/2}(\xx) =  \rho_f^{\frac{n+1}2} (s_l s_r)^{\frac12} s_f^{\frac{n+2}2} \Omega^{1/2}(\esxx),
$$
so that 
\begin{equation}\label{nubetamu}
\nu_0 = (\rho_l\rho_r)^{-\frac{n+1}2} \rho_f^{-(n+1)} (s_l s_r)^{-\frac{n+1}2} s_f^{-(n+1)} 
\beta^*\mu.
\end{equation}
Let $\beta_3$ be the blow down $\esxxx\to X^3$ and $\omega$ a smooth density on $\esxxx$.
Again, by computing codimensions we can see that
\begin{equation}\label{beta3mu}
\beta_3^* \mu_3^2 = (r_{LM} r_{MR} r_{LR})^{n+1}
r_{LMR}^{2n+2} (\gamma_L \gamma_M \gamma_R) (\gamma_{LM}\gamma_{MR}\gamma_{LR})^{n+2}
\gamma_{LMR}^{2n+4}\; \omega.
\end{equation}

Using the liftings (\ref{liftchart}) and combining the properties (\ref{mu3}), (\ref{nubetamu}), 
and (\ref{beta3mu}) we then obtain
$$
(\nu_0)_{LM} (\nu_0)_{MR} (\nu_0)_{LR} = r^{-(n+1)} \gamma^{-n} \omega,
$$
where $r$ is the product of all defining functions of the regular faces and $\gamma$ the same for the
crossover faces.  The composition formula (\ref{abnu}) then becomes
$$
(A\circ B) \rho^{-n} s^{-n+1}\; \omega_b
= \psi_{LR*} \Bigl[ A_{LM} B_{MR} r^{-n} \gamma^{-n+1} \omega_{b3} \Bigr],
$$
where $\rho = \rho_l\rho_r\rho_f$ and $s = s_l s_r s_f$.
The result follows because 
$$
\frac{r\gamma}{\psi_{LR}^*(\rho s)} = r_M \gamma_M.
$$ 
\end{proof}

Now let us prove the composition results needed to handle the Neumann series for $(I-F)^{-1}$.  Let $F'$ be
an asymptotic summation of the series $\sum_{l=1}^\infty F^l$ at the front crossover face.  Then $(\hDl) W
(I+F')$ will vanish to infinite order at both front faces.  The first step is to analyze $F^l$, for which we
need the following:
\begin{lemma}\label{FF}
Let $\calF = (\infty,\sigma_r,\infty,[\frac{n+1}2]_+,[\frac{n}2],[\frac12]_+)$.  Then
$\calA_\calF(\esxx; \shlf)$ is closed under composition.
\end{lemma}
\begin{proof}
With $A\nu_0, B\nu_0\in \calA_\calF(\esxx; \shlf)$, we will compute $(A\circ B)$ from the
formula (\ref{abomega}).  $A_{LM}$ and $B_{MR}$ are most succinctly described with a chart
of index sets akin to (\ref{liftchart}).  Fortunately we can consider the
regular and crossover faces separately.  At the regular faces we have:
\begin{equation}\label{abreg1}
\begin{tabular}{lccccccc}
\qquad & L\qquad& M\qquad& R\qquad& LM\qquad& MR\qquad&LR\qquad&LMR\\
$A_{LM}$&$\infty$&$\sigma_M$&--&$\infty$&$\sigma_R$&$\infty$&$\infty$\\
$B_{MR}$&--&$\infty$&$\sigma_R$&$\infty$&$\infty$&$\sigma_R$&$\infty$\\
\end{tabular}
\end{equation}
and at the crossover faces:
\begin{equation}\label{abcros1}
\begin{tabular}{lccccccc}
\qquad & L\qquad& M\qquad& R\qquad& LM\qquad& MR\qquad&LR\qquad&LMR\\
$A_{LM}$&$[\frac{n+1}2]_+$&$[\frac{n}2]$&--&$[\frac12]_+$&$[\frac{n}2]$&$[\frac{n+1}2]_+$&$[\frac12]_+$\\
$B_{MR}$&--&$[\frac{n+1}2]_+$&$[\frac{n}2]$&$[\frac{n+1}2]_+$&$[\frac12]_+$&$[\frac{n}2]$&$[\frac12]_+$\\
\end{tabular}
\end{equation}
Under $\psi_{LR}$ the faces which map to the interior are $\tilde M$ and $\tL_{M}$.  From (\ref{abcros1})
we see that the combination
$A_{LM} B_{MR}$ vanishes to infinite order at $\tilde M$, while at $\tL_{M}$ the index set is
$[n+\frac12]_+$, so the push-forward is well-defined.

Notice that $A_{LM} B_{MR}$ vanishes to infinite order 
at all of the regular faces of $\esxxx$ except $\tilde R$, where it is polyhomogeneous conormal 
with index $\sigma_R$.  Thus $\rho_r^{-\tsr} (A\circ B)\omega_b$ will be the pushforward of a
b-density with index set $\{(k,k); k\in\bbN_0\}$ at $\tilde R$ and $\infty$ at all the other regular
faces.  By Proposition \ref{pushfor} we conclude that $(A\circ B)$ has orders
$\infty,\sigma_r,\infty$ at the left, right, and front regular faces.

As for the crossover faces, by (\ref{abcros1}) and Proposition \ref{pushfor} 
the index sets of $(A\circ B)$ are $[\frac{n+1}2]_+ \xU [\frac{n}2+1]_+$ at the left crossover,
$[\frac{n}2] \xU [\frac{n+1}2]_+$ at the right crossover, and $[1]_+ \xU [n+\frac12]_+$ at the front
crossover face.  It is a trivial observation that for $a<b$
$$
[a] \xU [b]_+ = [a],\qquad [a]_+ \xU [b]_+ = [a]_+,
$$ 
and the result follows.
\end{proof}

Applying this Lemma to the powers of $F$, we may conclude that
$$
F' \in \calA_{\calF}(\esxx;\shlf); \qquad \calF =
(\infty,\sigma_r,\infty,[\tfrac{n+1}2]_+,[\tfrac{n}2],[\tfrac12]_+).
$$
Actually, this doesn't quite follow directly from Lemma \ref{FF}.  We need to know in addition that
at each order there is a bound on the highest logarithmic power in $F^l$ at the front crossover face which
is  uniform for $l\in \bbN_0$.    This can be deduced from the fact that
the index set of $F^l$ at the front crossover face is really $[\frac{l}2]_+$, as one can see from the
proof of  Lemma \ref{FF}.

The improved parametrix is $M = W (I+F')$, which by the construction of $F'$ satisfies
\begin{equation}\label{pareq}
I - (\hDl) M = E \in \calA_\calE(\esxx; \shlf),
\end{equation}
where $\calE = (\infty,\sigma_r,\infty, [\frac{n+1}2]_+, [\frac{n}2], \infty)$.
\begin{proposition}\label{compprop}
For $\calM$, $\calF$ as above we have
\begin{equation}\label{compos}\begin{split}
\psop^{*}(X; \shlf) \circ \calA_\calF(\esxx; \shlf) &\subset \calA_\calF(\esxx; \shlf)\\
\calA_\calM(\esxx; \shlf) \circ \calA_\calF(\esxx; \shlf) &\subset \calA_{\calM\circ\calF}(\esxx; \shlf),
\end{split}\end{equation}
where $(\calM\circ\calF) = (\sigma_l,\sigma_r,2\sigma_f,[\frac{n}2],[\frac{n}2],[\frac12]_+)$.  
\end{proposition}
\begin{proof}
Let us start with the second formula in (\ref{compos}).  Suppose that 
$$
A \nu_0 \in \calA_\calM(\esxx; \shlf), \qquad B \nu_0 \in  \calA_\calF(\esxx; \shlf).
$$  
As in the proof of Lemma \ref{FF}, we examine the behavior of the lifts to $\esxxx$. 
At the regular faces we have:
\begin{equation}\label{abreg2}
\begin{tabular}{lccccccc}
\qquad & L\qquad& M\qquad& R\qquad& LM\qquad& MR\qquad&LR\qquad&LMR\\
$A_{LM}$&$\sigma_L$&$\sigma_M$&--&$0$&$\sigma_R$&$\sigma_R$&$0$\\
$B_{MR}$&--&$\infty$&$\sigma_R$&$\infty$&$\infty$&$\sigma_R$&$\infty$\\
\end{tabular}
\end{equation}
In order to apply Proposition \ref{pushfor} we consider $x^{-\tsl} {x'}^{-\tsr} (A\circ B)$, which is the
push-forward of a b-density on $\esxxx$ with index set $\{(k,k);\;k\in\bbN_0\}$ at 
$\tilde L$, $\tilde R$, and $\tilde S_{LR}$, and $\infty$ otherwise.  
Then pushforward gives indices $(\sigma_l, \sigma_r, 2\sigma_f)$ for $A\circ B$ at the regular faces.

The index chart for the crossover faces:
\begin{equation}\label{abcros2}
\begin{tabular}{lccccccc}
\qquad & L\qquad& M\qquad& R\qquad& LM\qquad& MR\qquad&LR\qquad&LMR\\
$A_{LM}$&$[\frac{n}2]$&$[\frac{n}2]$&--&$0$&$[\frac{n}2]$&$[\frac{n}2]$&$0$\\
$B_{MR}$&--&$[\frac{n+1}2]_+$&$[\frac{n}2]$&$[\frac{n+1}2]_+$&$[\frac12]_+$&$[\frac{n}2]$&$[\frac12]_+$\\
\end{tabular}
\end{equation}
We conclude that the crossover index sets for $A\circ B$ are $([\frac{n}2],[\frac{n}2],[\frac12]_+)$.

The first formula in (\ref{compos}) is proven in exactly the same way.  If $B$ is as above, but
$A\nu_0 \in \psop^{*}(X; \shlf)$, 
then the product $A_{LM} B_{MR}$ will have an interior singularity at the lift of the
diagonal through $\psi_{LM}$.  One can check that $\psi_{LR*}$ annihilates this singularity by
standard wave-front set arguments (see \cite{EMM} for an explicit discussion of this).  Then the proof of
push-forward to $\calA_\calM$ is exactly as above.
\end{proof}

\begin{corollary}\label{mstruct}
The operator $M = W(I+F')$ satisfies
$$
M \in \psop^{-2}(X; \shlf) + \calA_{\calM}(\esxx; \shlf).
$$
\end{corollary}
\begin{proof}
Let $W = W_1 + W_2$ as an element of $\psop^{-2}(X; \shlf) + \calA_{\calM}(\esxx; \shlf)$.
From (\ref{compos}) we immediately see that 
\begin{equation*}\begin{split}
W_2 + W_1 F' &\in \calA_{\calM}(\esxx; \shlf),\\
W_2 F' &\in \calA_{(\calM\circ\calF)}(\esxx; \shlf).
\end{split}\end{equation*}
Now consider the equation (\ref{pareq}).  By construction $I - (\hDl) W_1 = Q$, which is smooth
up to the front face, and since $\Delta$ lifts to $\Diff_b^*(\esxx)$ we conclude that
$$
(\hDl) (W_2 F') \in \calA_{\calM}(\esxx; \shlf).
$$
This means that $\hDl$ must annihilate all of the terms in the asymptotic expansion
of $(W_2 F')$ at the front regular face, since an expansion containing $\rho_f^{2\sigma_f}$
is not allowed in $\calA_{\calM}$.  We can argue term by term using the normal operator
to show that these coefficients are zero.  For example, at leading order we have
\begin{equation}\label{nwf}
N_q(\hDl) N_q({x'}^{-2\sigma_f}W_2 F') = 0,
\end{equation}
for all $q\in \bX\bck \Lambda$.
From the definition of $\calA_{\calM\circ \calF}(\esxx; \shlf)$ we see that 
$$
N_q({x'}^{-2\sigma_f}W_2 F') \in \calA_{\sigma(q), \sigma(q)}(F_q; \shlf). 
$$
Applying Proposition \ref{unique} (in the model case with constant indicial root) to (\ref{nwf}), 
we conclude that $N_q({x'}^{-2\sigma_f}W_2 F') =0$.  Using this argument inductively, we
conclude that $W_2 F'$ vanishes to infinite order at the front face, which puts it in 
$\calA_\calM(\esxx; \shlf)$.
\end{proof}

Next we show that the structure of the error term $E$ allows it to be realized as a conormal function
on $\xlxl$ rather than $\esxx$.  Because we can realize $\esxx$ as $[\xx;
\Lambda_l; \Lambda_r; \Lambda_f; S]$, by various applications of Lemma \ref{buorder}, there is a
well-defined blow-down $\phi: \esxx\to \xlxl$.  This is a b-map, but not a b-fibration
($\tL_f$ is mapped into a corner, for example).  
\begin{lemma}
Suppose as above that
$$
E \in \calA_\calE(\esxx; \shlf), \qquad \calE = (\infty,\sigma_r,\infty, [\tfrac{n+1}2]_+,
[\tfrac{n}2], \infty).
$$  
Then if we lift $E$ from the interior $(\xx)^\circ$ to $\xlxl$ we have
$$
E \in \calA_\calH(\xlxl; \shlf), \qquad \calH = (\infty, \sigma_r, [\tfrac{n+1}2]_+, [\tfrac{n}2])
$$
\end{lemma}
\begin{proof}
First note that $\calV_b(\xlxl)$ lifts under $\phi$ to $(\rho_f s_f)^{-1} \calV_b(\esxx)$,
as may be easily checked in local coordinates.
Observe that 
$$
x^{-\frac{n+1}2} {x'}^{-\tsr} E\in \calA^{(\infty,0,\infty,0,0,\infty)-}(\esxx).
$$ 
Because of the infinite order of vanishing at the two front faces, we deduce that
$$
x^{-\frac{n+1}2} {x'}^{-\tsr} E\in \calA^{(\infty,0,0,0)-}(\xlxl)
$$

Let $V_l, V_r, W_l, W_r$ be radial vector fields at the left/right regular and left/right crossover
boundary faces of $\xlxl$, respectively.  We claim the lift by $\phi$ of
$V_l$ is a radial vector field $\tilde V_l$ at the left regular face, possibly plus a term in
$\rho_l s_f^{-1} \calV_b(\esxx)$.  Similarly, $W_l$ lifts
to a radial vector field $\tilde W_l$ at the left crossover face, possibly plus a term in
$s_l \rho_f^{-1} \calV_b(\esxx)$.   The behavior on the right is analogous.
These facts are best checked in local coordinates.  On $\xlxl$ we can use coordinates
$(u,\theta,\bar y,u',\theta',\bar y')$, such that $x = u\cos \theta$, $y_n = u\sin\theta$,
$y = (\bar y,y_n)$, and similarly on the right.  For $V_l$ and $W_l$ we can take $(\cos\theta)\del_\theta$
and $u\del_u$.

Near $\tL_l$ and $\tL_r$ we ignore the blow-up of $S$, since $\tilde S$ does not intersect these.
So we can use the coordinates $(R,\eta,\eta',\omega,\theta,\theta',\bar y)$, where 
$R = \sqrt{r^2+{r'}^2+(\bar y - \bar y')^2}$ and $(\eta,\eta',\omega) = (r,r',\bar y-\bar y')/R$.
Then $\tilde W_l = (\cos\theta)\del_\theta$ again, which is indeed a radial vector field at
the left regular face $\{\cos\theta = 0\}$. The lift of $V_l$ is
$$
\tilde V_l = \eta\del_\eta + \eta^2 \bigl[R\del_R - \eta\del_\eta - \eta'\del_{\eta'} - 
\omega\del_\omega\bigr],
$$
which is a radial vector field at the left crossover face $\{\eta = 0\}$.  
Similarly, with a different coordinate system we could check this in a neighborhood of $S$ (then we could
ignore the blow-ups of $\Lambda_l$ and $\Lambda_r$).

Of course, the radial vector fields are not uniquely defined.  Another choice of $V_l$ would
differ by a field in $t_l \calV_b(\xlxl)$, where $t_l$ is a defining function for the left regular face of
$\xlxl$. Since $\phi^*t_l = \rho_l \rho_f$, $t_l \calV_b(\xlxl)$ lifts to $\rho_l s_f^{-1} \calV_b(\esxx)$,
which agrees with the correction term given above.

Now we can show $E \in \calA_\calH(\xlxl; \shlf)$ by pulling these radial vector fields up to
$\esxx$.  (Once again the infinite order of vanishing at the front faces is crucial.)
For example, $E \in \calA_\calE(\esxx; \shlf)$ entails the estimate
$$
\tilde V_r (x^{-\frac{n+1}2} {x'}^{-\tsr} E) \in \calA^{(\infty,1,\infty,0,0,\infty)-}(\esxx),
$$
and this would be unaffected if we add to $\tilde V_r$ a term in $\rho_r s_f^{-1} \calV_b(\esxx)$.
Consequently, it implies the estimate
$$
V_r (x^{-\frac{n+1}2} {x'}^{-\tsr} E) \in \calA^{(\infty,1,0,0)-}(\xlxl).
$$
In a similar way all of the estimates needed to show $E \in \calA_\calH(\xlxl; \shlf)$ can
be deduced.
\end{proof}

We consider next the composition of a kernel living on 
$\esxx$ with a kernel on $\xlxl$. For this purpose we can use the slightly simpler space
$$
(\esxx)\times \Xl = [X^3; \Lambda_{LM}; S_{LM}; S_{LR};  \Lambda_{L}; \Lambda_{M}; \Lambda_{R}]
$$ 
instead of $\esxxx$ for the triple product.  We'll continue to denote the projections onto 
pairs by $\psi_{LM},\psi_{MR},\psi_{LR}$, where $\psi_{LM}$ projects to $\esxx$ and
$\psi_{MR}$ and $\psi_{LR}$ are blow-downs to $(\Xl)^3$ followed by projection onto $\xlxl$.
The chart (\ref{liftchart}) of lifts of defining functions is abbreviated to
\begin{equation}\label{liftchart2}
\begin{tabular}{lcccc}
\qquad & L\qquad& M\qquad& R\qquad& LM\\
$\psi_{LM}^*$&l&r&--&f\\
$\psi_{MR}^*$&--&l&r&l\\
$\psi_{LR}^*$&l&--&r&l\\
\end{tabular}
\end{equation}
Let $\nu_0$ be a section of $\shlf(\esxx)$ and $\mu_0$ a section of 
$\shlf(\xlxl)$.  The push-forward formula for composition of $A \nu_0$ with $B\mu_0$ is
$$
(A\circ B)\mu_0 = \psi_{LR*} \Bigl[ A_{LM} B_{MR} (\nu_0)_{LM} (\mu_0)_{MR} \Bigr].
$$
Computing as in the proof of Lemma \ref{compdef}, we derive
\begin{equation}\label{respush}
(A\circ B)\omega_b = \psi_{LR*} \Bigl[ A_{LM} B_{MR} \rho_M^{-n} \gamma_M^{-n+1} \omega_{b3} \Bigr],
\end{equation}

\begin{lemma}\label{residlem}
With 
$$
\calM = (\sigma_l, \sigma_r, 0, [\tfrac{n}2], [\tfrac{n}2], [0]),\qquad
\calH = (\infty, \sigma_r, [\tfrac{n+1}2]_+, [\tfrac{n}2]),
$$ 
as above, we have the compositions
\begin{equation*}\begin{split}
\calA_\calM(\esxx; \shlf) &\circ \calA_\calH(\xlxl; \shlf) \subset \calA_\calI(\xlxl; \shlf) \\
\psop^{-2}(X; \shlf) &\circ \calA_\calH(\xlxl; \shlf) \subset \calA_\calH(\xlxl; \shlf) \\
\end{split}\end{equation*}
where $\calI = (\sigma_l, \sigma_r, [\frac{n}2], [\frac{n}2])$.
\end{lemma}
\begin{proof}
Let $A \nu_0 \in \calA_\calM(\esxx; \shlf)$ and $B\mu_0 \in \calA_\calH(\xlxl; \shlf)$.
As in the proof of Lemma \ref{FF}, we use a chart to describe the boundary behavior 
of $A_{LM} B_{MR}$.  For the regular faces we have:
$$
\begin{tabular}{lcccc}
\qquad & L\qquad& M\qquad& R\qquad& LM\\
$\psi_{LM}^*$&$\sigma_L$&$\sigma_M$&--&$\sigma_{LM}$\\
$\psi_{MR}^*$&--&$\infty$&$\sigma_R$&$\infty$\\
\end{tabular}
$$
and at the crossover faces:
$$
\begin{tabular}{lcccc}
\qquad & L\qquad& M\qquad& R\qquad& LM\\
$\psi_{LM}^*$&$[\tfrac{n}2]$&$[\tfrac{n}2]$&--&$[\tfrac{n}2]$\\
$\psi_{MR}^*$&--&$[\tfrac{n+1}2]_+$&$[\tfrac{n}2]$&$[\tfrac{n+1}2]_+$\\
\end{tabular}
$$
The behavior at the middle faces, together with (\ref{respush}), shows that
the push-forward is well defined.  And the form of $\calI$ is deduced from
(\ref{liftchart2}) and Proposition \ref{pushfor}.

The second formula is handled in the same way.  As in Proposition \ref{compprop} the push-forward
annihilates the interior singularity.
\end{proof}

As a final topic of this subsection, we establish a mapping property that will be used in the
characterization of generalized eigenfunctions.  
\begin{lemma}\label{xlmap}
With,
$$
\calM = (\sigma_l, \sigma_r, 0, [\tfrac{n}2], [\tfrac{n}2], [0]), \qquad
\calI = (\sigma_l, \sigma_r, [\tfrac{n}2], [\tfrac{n}2]),
$$
we have the mapping properties:
\begin{equation*}\begin{split}
\psop^{*}(X; \shlf) \circ \calA_{\infty, [\frac{n+1}2]_+}(\Xl) &\subset \calA_{\infty,
[\frac{n+1}2]_+}(\Xl)\\
\calA_{\calM}(\esxx;\shlf) \circ \calA_{\infty, [\frac{n+1}2]_+}(\Xl) &\subset \calA_{\sigma,
[\frac{n}2]}(\Xl)\\
\calA_{\calI}(\xlxl;\shlf) \circ \calA_{\infty, [\frac{n+1}2]_+}(\Xl) &\subset \calA_{\sigma,
[\frac{n}2]}(\Xl)\\
\end{split}\end{equation*}
\end{lemma}
\begin{proof}
We will prove only the second of the three formulas.  The first follows by similar argument because of the
cancelation of the interior singularity as in Proposition \ref{compprop}.  The third formula follows
by an analogous but even simpler argument on $\xlxl$. 

Let $(A\nu_0) \in \calA_{\calM}(\esxx;\shlf)$, where $\nu_0$ is a section of $\shlf(\esxx)$.
As noted in the remarks following (\ref{esxxdef}), there are b-fibrations
$$
\beta_L, \beta_R: \esxx\to\Xl
$$
projecting through the left and right factors.
The action of the operator with kernel $A \nu_0$ on $u \gamma_0 \in\calA_{\infty, [\frac{n+1}2]_+}(\Xl)$,
where $\gamma_0\in \cinf(\Xl;\shlf)$, is given by 
$$
(A\circ u) \gamma_0 = \beta_{L*} \Bigl[ A \nu_0\cdot \beta_R^*(u \gamma_0)\Bigr]
$$
Multiplying both sides by $\gamma_0$, we can recast this in a form suitable for 
push-forward:
$$
(A\circ u) \mu_b = \beta_{L*} \Bigl[ A u_R\cdot \rho_r^{-n} s_r^{-n+1}\omega_b\Bigr],
$$ 
where $u_R$ is shorthand for $\beta_R^* u$, $\omega_b\in \cinf(\esxx,\Omega_b)$, and $\mu_b \in
\cinf(\Xl;\Omega_b)$.  The combination $A u_R$ has index family $(\sigma_l, \infty,\infty,[\frac{n}2],
[\frac{n+1}2]_+, [\frac{n+1}2]_+)$, so the push-forward is well-defined.  Multiplying
$(A\circ u)$ by $x^{-\tsigma}$ and applying the push-forward result Proposition \ref{pushfor}
tells us that $x^{-\tsigma} (A\circ u)$ has index family $(0, [0]\xU [\frac{1}2]_+)$, and the result
follows.
\end{proof}

\subsection{Resolvent kernel}
Consider the operator $\Dl$ restricted to a domain defined by the radiation condition:
$$
\calD_\lambda = \{u\in x^{-\delta}L^2(X;\shlf): \; (x\del_x - \sigma)u \in x^\delta L^2(X;\shlf)\} \cap
H^2(X;\shlf),
$$
for some small $\delta>0$.
Thus
$$
\hDl: \calD_\lambda \to x^\delta L^2(X;\shlf).
$$
The limiting absorption principle shows the existence of an inverse
$$
R_\zeta:  x^\delta L^2(X;\shlf) \to \calD_\lambda,
$$
but we do not know the structure of its kernel.

On the other hand, in \S\ref{compss} we obtained the parametrix $M$ such that
$$
(\hDl) M = 1-E \text{ on }x^\delta L^2(X;\shlf),
$$
where we do know the structure of $M$ and $E$.
For $u\in \calD_\lambda$, integration by parts is justified so that we can take the transpose
$$
M^t (\hDl) = 1-E^t,
$$
as a relation on $\calD_\lambda$.  With these relations we can write
\begin{equation*}\begin{split}
M &= R_\zeta(\Dl)M = R_\zeta-R_\zeta E\\
M^t &= M^t(\Dl)R_\zeta = R_\zeta -E^tR_\zeta
\end{split}\end{equation*}
And so we can deduce
\begin{equation}\label{rmre}
R_\zeta = M+R_\zeta E = M + M^t E + E^t R_\zeta E.
\end{equation}

By Corollary \ref{mstruct}
$$
M^t \in \psop^{-2}(X; \shlf) + \calA_{\calM}(\esxx; \shlf).
$$
(Taking the transpose of an index family interchanges left and right index sets, but $\calM$
is symmetric.)  Then by Lemma \ref{residlem} we have
\begin{equation}\label{mte}
M^t E \in \calA_\calI(\xlxl; \shlf)(X; \shlf),\qquad 
\calI = (\sigma_l, \sigma_r, [\tfrac{n}2], [\tfrac{n}2]).
\end{equation}
 
To deal with $E^t R_\zeta E$, we can at least use the fact that $R_\zeta$ is bounded as a map $x^\delta L^2
\to x^{-\delta} L^2$ for small $\delta>0$.  

\begin{lemma}\label{resbndd}
\begin{equation*}\begin{split}
&\calA_{\calH^t}(\xlxl; \shlf) \circ \calL(x^\delta L^2, x^{-\delta} L^2) \circ 
\calA_\calH(\xlxl; \shlf)\\ 
&\qquad\subset \calA_\calI(\xlxl; \shlf)(X; \shlf),
\end{split}\end{equation*}
where $\calI = (\sigma_l, \sigma_r, [\frac{n}2], [\frac{n}2])$.
\end{lemma}
\begin{proof}
Since the index family $\calH$ implies $x^\delta L^2$ decay at the left face, the composition
is represented by a convergent integral.  The necessary estimates can be obtained by differentiating under
the integral using dominated convergence.
\end{proof}

With (\ref{mte}) and Lemma \ref{resbndd} we deduce immediately from (\ref{rmre}) the main result of this
section:
\begin{theorem}\label{resdirect}
For $\zeta$ such that $\lambda = \alpha_0^2\zeta(n-\zeta) \in (\frac{\alpha_0^2n^2}4, 
\frac{\alpha_1^2n^2}4)$ and satisfying the generic assumption (\ref{crassume}),
the kernel of the resolvent $R_\zeta = [\hDl]^{-1}$ defined by the limiting absorption principle has the
structure
$$
R_\zeta \in \psop^{-2}(X; \shlf) + \calA_{\calM}(\esxx; \shlf) + \calA_{\calI}(\xlxl; \shlf),
$$
with $\calM = (\sigma_l,\sigma_r,0,[\tfrac{n}2],[\tfrac{n}2],0)$ and 
$\calI = (\sigma_l,\sigma_r,[\tfrac{n}2],[\tfrac{n}2])$. 
\end{theorem}

With this knowledge of the structure of the resolvent kernel, we can give a more refined version 
of Proposition \ref{nuparam}, describing the generalized eigenfunctions.  
\begin{proposition}
For $\zeta, \lambda$ as in Theorem \ref{resdirect}, given $f \in \cinf(\bX)$ we can solve $(\Dl) u =0$ for 
$$
u \in \calA_{n-\sigma,[\tfrac{n}2]}(\Xl) + \calA_{\sigma,[\tfrac{n}2]}(\Xl),
$$
such that
$$
u \sim x^{n-\sigma} f + x^{\sigma} f' \text{ near }W_\lambda,
$$
with $f' \in \calA_{[\tfrac{n}2]}(\overline{W_\lambda})$.  This $u$ is uniquely determined by
$f|_{\overline{W_\lambda}}$
\end{proposition}
\begin{proof}
Using Proposition \ref{xlsolve} we solve for $u_1 \in \calA_{n-\sigma,[\tfrac{n}2]}(\Xl)$ and asymptotic to
$x^{n-\sigma}f$ near $W_\lambda$, such that $(\Dl) u_1 = \phi \in \calA_{\infty, [\tfrac{n+1}2]_+}(\Xl)$.  
Then let $u_2 = -R_\zeta\phi$, which is in $\calA_{\sigma,[\tfrac{n}2]}(\Xl)$ by Lemma \ref{xlmap}.
Thus $u = u_1 + u_2$ has the stated properties.
\end{proof}

In particular, under these assumptions the scattering matrix defined in Theorem \ref{nusmatrix}
extends to a map
$$
S_\lambda: \cinf(\overline{W_\lambda}) \to \calA_{[\tfrac{n}2]}(\overline{W_\lambda}).
$$

\end{document}